\documentclass[graybox]{svmult}
\usepackage{mathptmx}       
\usepackage{helvet}         
\usepackage{courier}        
\usepackage{makeidx}         
\usepackage{graphicx}        
\usepackage{multicol}        
\usepackage[bottom]{footmisc}
\usepackage{amsmath,amssymb,amsfonts}        
\usepackage{url}
\usepackage{pgf}
\usepackage{color}
\usepackage{hyperref}
\usepackage{accents}
\usepackage{array}
\makeindex             

\def\term{\textbf}

\renewcommand{\D}{\mathbb{D}}

\newcommand{\F}{\mathcal{F}}

\newcommand{\fpqr}{\phi^t_{pqr}}
\newcommand{\fsus}{\phi_{\mathrm{susp}}}

\newcommand{\Gpqr}{\Gamma_{pqr}}

\newcommand{\Hy}{\mathbb{H}^2}

\newcommand{\R}{\mathbb{R}}

\newcommand{\Sph}{\mathbb{S}}
\renewcommand{\S}{\mathcal{O}}
\newcommand{\Spqr}{\S_{pqr}}
\newcommand{\SLZ}{\mathrm{SL}_2(\Z)}

\newcommand{\tr}{\mathrm{tr}}

\newcommand{\U}{\mathrm{T}^1}

\newcommand{\Z}{\mathbb{Z}}

\newcommand{\catbat}{M~~~
	\begin{picture}(4,0)
	\put(-8,-2.3){\includegraphics[height=2.7mm]{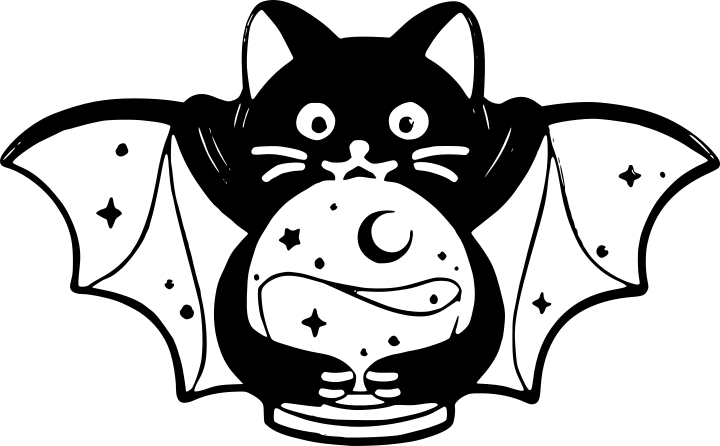}}
	\end{picture}
}

\begin{document}
\title*{The cat-bat map, the figure-eight knot, and the five orbifolds}
\author{Pierre Dehornoy}
\institute{Aix-Marseille Université, CNRS, I2M, Marseille, France, \email{pierre.dehornoy@univ-amu.fr}}
%
%
\maketitle

\begin{center}
\includegraphics[width=.8\textwidth]{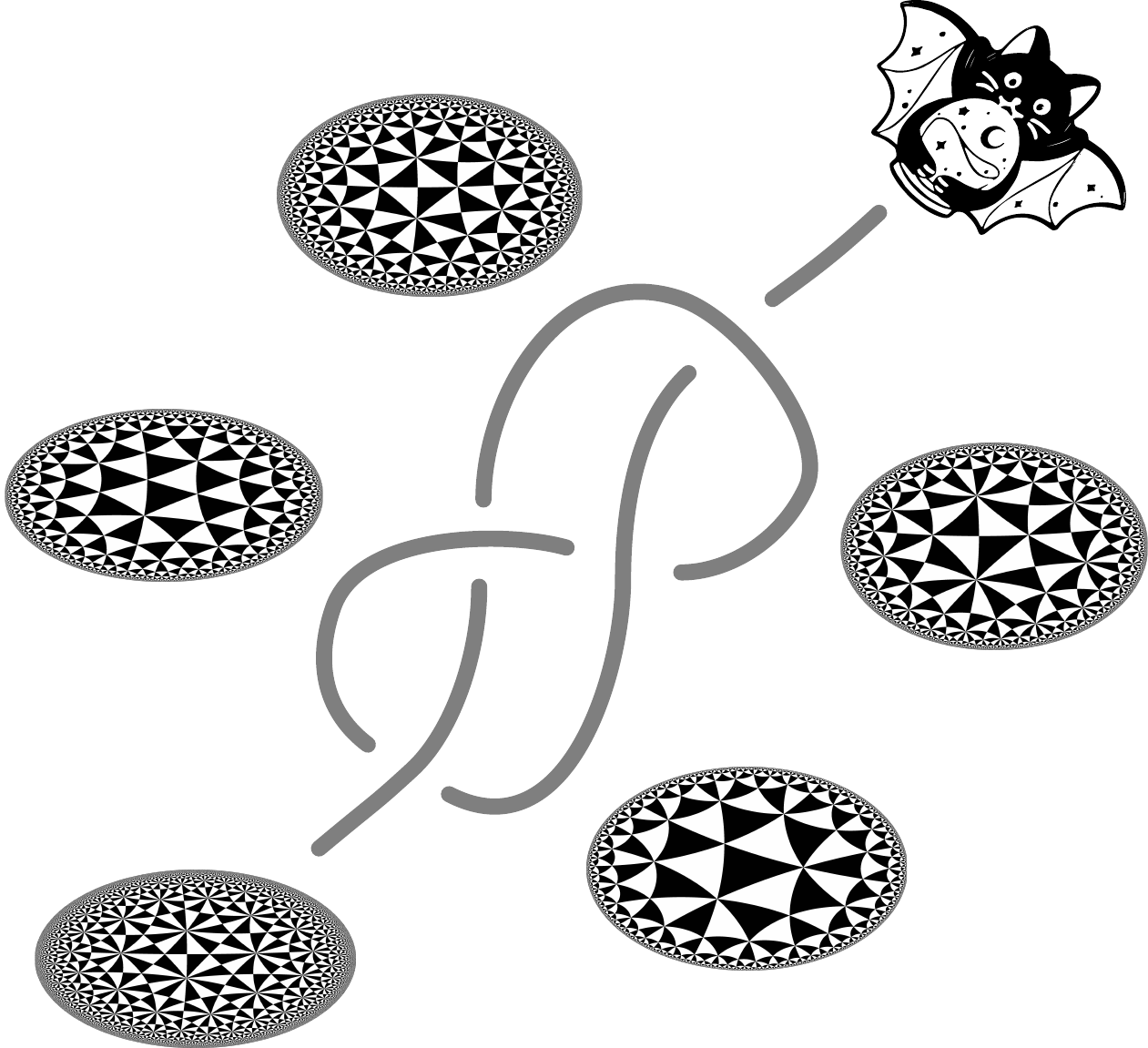}
\end{center}

One of the most famous ways to construct examples hyperbolic 3-manifolds is by Dehn filling the complement of the {figure-eight} knot in the 3-sphere. 
Indeed William Thurston showed that all but 11 fillings give rise to hyperbolic 3-manifolds~\cite[Chap 4]{ThurstonNotes}. 
For these exceptions, one may wonder which manifolds are obtained. 
Since the figure-eight knot is fibered with fiber a once-punctured torus, and since the associated monodromy map is the so-called\footnote{Boris Hasselblatt noticed, in Chapter 3 ``Hyperbolic dynamical systems'' of the Handbook of Dynamical Systems 1A, that it is not clear whether Arnold's picture of a cat is the oldest picture of that map. Another possibility is Avez' bat. Another possibility is to speak of the auThomorphism of the torus, since René Thom seems to have noticed his properties before Arnold and Avez. We stick to the cat-bat.} \term{cat-bat map}~$(\begin{smallmatrix}2&1\\1&1\end{smallmatrix}):\R^2/\Z^2\to \R^2/\Z^2$, the 0-surgery yields the mapping torus of the cat-bat map. 
If you ask a software like Regina\footnote{This is what Andy Hammerlindl did, and that originates this note.}, she will tell you that the $\pm1$-, $\pm1/2$-, and $\pm1/3$-surgeries yield the unit tangent bundle to the some hyperbolic 2-orbifolds and that the $\pm1/4$-surgeries yield a graph-manifold.

One can rephrase the previous observation as follows. 
Denote by $\catbat$ the mapping torus of the cat-bat map, equipped with the suspension flow~$\fsus^t$. 
The cat-bat map has a unique fixed point, coming from the point~$(0,0)$ in~$\R^2$. 
Call $\gamma_1$ the orbit of~$\fsus^t$ corresponding to this unique fixed point.
Write $\catbat(\gamma_1, b/a)$ for the 3-manifold obtained by Dehn surgery\footnote{We use Rolfsen's convention for the surgery coefficient~$b/a$, meaning the meridian of the new plain torus is homologous to~$a\ell+b m$, where $\ell$ is a longitude --given here by following~$\gamma_1$ along the invariant stable direction-- and $m$ a meridian. In particular the denominator of the slope is the number of intersection points between the old and the new meridians.} on~$\gamma_1$ with slope $b/a$. 
Then $\catbat(\gamma_1, b/a)$ is hyperbolic, except when $b/a=0,\pm1,\pm\frac12,\pm\frac13,\pm\frac14$, and for $b/a=\pm1,\pm\frac12,\pm\frac13$ the manifold is the unit tangent bundle to a hyperbolic orbifold. 

One can further use Regina to explore surgeries on other orbits of the suspension flow of the cat-bat map. 
For example, the cat-bat map has two orbits of period 2 (corresponding to the points $(3/5, 1/5)$ and $(1/5, 2/5)$ in~$\R^2/\Z^2$). 
Denote by $\gamma_2$ the orbit of the suspension flow going through any of these two order 2-points, then $\catbat(\gamma_2, b/a)$ is unit tangent bundle to a hyperbolic orbifold for $b/a=\pm1$ and $\pm\frac12$. 

Our goal is to picture these five exceptional surgeries as explicitly as possible. 
Denote by~$\S_{pqr}$ the hyperbolic triangular orbifold which is a sphere with three cone points of angles $2\pi/p, 2\pi/q, 2\pi/r$ and by $\U\S_{pqr}$ its unit tangent bundle which is a 3-manifold, see Section~\ref{S:Orbifolds}.

\begin{theorem}\label{T:main}
The manifold $\catbat(\gamma_1, b/a)$ is 
\begin{itemize}
\item $\U\S_{237}$ when $b/a=\pm1$,
\item $\U\S_{245}$ when $b/a=\pm1/2$,
\item $\U\S_{334}$ when $b/a=\pm1/3$;
\end{itemize}
The manifold $\catbat(\gamma_2, p/q)$ is 
\begin{itemize}
\item $\U\S_{246}$ when $b/a=\pm1$,
\item $\U\S_{344}$ when $b/a=\pm1/2$.
\end{itemize}
\end{theorem}

As explained above, a software like Regina can easily prove Theorem~\ref{T:main}. 
Our point is to explain how one can {\it see} these surgeries. 
We do it in the reverse direction as we now explain: like any unit tangent bundle, each of the five manifolds~$\U\S_{pqr}$ we are interested in carries a geodesic flow, that we denote by~$(\fpqr)_{t\in\R}$. 
Since the orbifolds we study are good and hyperbolic ({\it i.e.} quotients of the hyperbolic plane), the geodesic flow is of Anosov type.
For each of these geodesic flows we find a genus-1-Birkhoff section, that is, a surface with boundary~$S$ which is a torus, which is bounded by orbits of~$\fpqr$, whose interior is transverse to~$\fpqr$, and which intersects all orbits of~$\fpqr$ in bounded time. 
In that case there is a well-defined first-return map on the interior of~$S$ along~$\fpqr$. 
When contracting each boundary component of~$S$ into a point, we obtain a well-defined map~$f_{\bar S}:\bar S\to\bar S$. 
The facts that $\fpqr$ is of Anosov type with coorientable invariant foliations and that $S$ is of genus 1 together imply that $f_{\bar S}$ is of Anosov type. 
Therefore it is conjugated to a linear map of the torus. 
A careful analysis then shows that $f_{\bar S}$ has exactly one fixed point. 
It turns out that, up to conjugacy, there is only one such map of the torus, namely the cat-bat map. 
This implies that the surgery on~$\U\S_{pqr}$ along the orbits which are the boundary of the surface~$S$ with an appropriate direction (which depends on the shape of~$S$ around~$\partial S$) yields the manifold~$M$. 
Thanks to the above argument, the technical property we explain and which implies Theorem~\ref{T:main} is

\begin{proposition}
\begin{enumerate}
\item For $(p, q, r)=(2,3,7)$ $(${\it resp.} $(2,4,5), (3,3,4)$$)$,  the geodesic flow $\fpqr$ on~$\U\Spqr$ admits a genus-one Birkhoff section~$S_{pqr}$ with one boundary component and boundary direction $(1,1)$ $(${\it resp.} $(2, 1), (3,1))$, such that the induced first-return map on~$\bar S_{pqr}$ has exactly one fixed point.
\item For $(p, q, r)=(2,4,6)$ $(${\it resp.} $(3,4,4))$, the geodesic flow~$\fpqr$ on~$\U\Spqr$ admits a genus-one Birkhoff section~$S_{pqr}$ with two boundary components on the same periodic orbit and boundary direction $(2,2)$ $(${\it resp.} $(4,2))$, such that the induced first-return map on~$\bar S_{pqr}$ has exactly one fixed point.
\end{enumerate}
\end{proposition}

The strategy we follow is not really new. 
It goes back to Fried, Ghys, Hashiguchi, Brunella who already observed and used the fact that some geodesic flows admit surgeries that yield suspensions of toral automorphism~\cite{Fried, Ghys, Hashiguchi, Brunella}. 
The first case of Theorem~\ref{T:main} was explicitly stated about 10 years ago~\cite{D:ETDS} (and probably known to some experts even earlier).
The third case was suggested by Pierre Will, see~\cite{Will} for a connection with complex hyperbolic geometry. 
It was explained using different and less direct techniques (computations of linking numbers) in a recent note~\cite{D:TSG}. 
The second, fourth and fifth cases were suggested by Andy Hammerlindl. 
These cases use new variants of the constructions of~\cite{D:ETDS}.

\bigskip

\noindent {\bf Acknowledgments.} 
This note emerged from discussions with Andy Hammerlindl during the Matrix event ``Dynamics, Foliation, and Geometry III'' at the mathematical research institute Matrix.
I thank Andy for these motivating discussions, as well as Andy and Jessica Purcell for organising this very nice event, and Matrix for hosting it. 
The 334- and 344-cases were also discussed with Neige Paulet, I thank her for her remarks. 
I also thank the referee for several suggestions that improve the readability of the paper. 
I acknowledge support from ANR Gromeov ANR-19-CE40-0007, and I thank Michele Triestino for leading~it.


\section{Hyperbolic orbifolds, geodesic flows, and Birkhoff sections}
\label{S:Prelim}

\subsection{Hyperbolic triangular orbifolds and their unit tangent bundles}\label{S:Orbifolds}
Let $p, q, r$ be three integers satisfying $1/p+1/q+1/r<1$. 
In the hyperbolic plane~$\Hy$, there is a triangle~$PQR$ whose angles are $\pi/p, \pi/q, \pi/r$, respectively.
Consider the group~$\Gpqr$ of hyperbolic isometries generated by the rotations of respective angles $2\pi/p, 2\pi/q, 2\pi/r$ around $P, Q, R$. 
The action of~$\Gpqr$ on~$\Hy$ is proper and cocompact. 
A fundamental domain is obtained by taking the union of~$PQR$ with its mirror image by the symmetry around any of its sides. 

\begin{center}
\includegraphics[width=.3\textwidth]{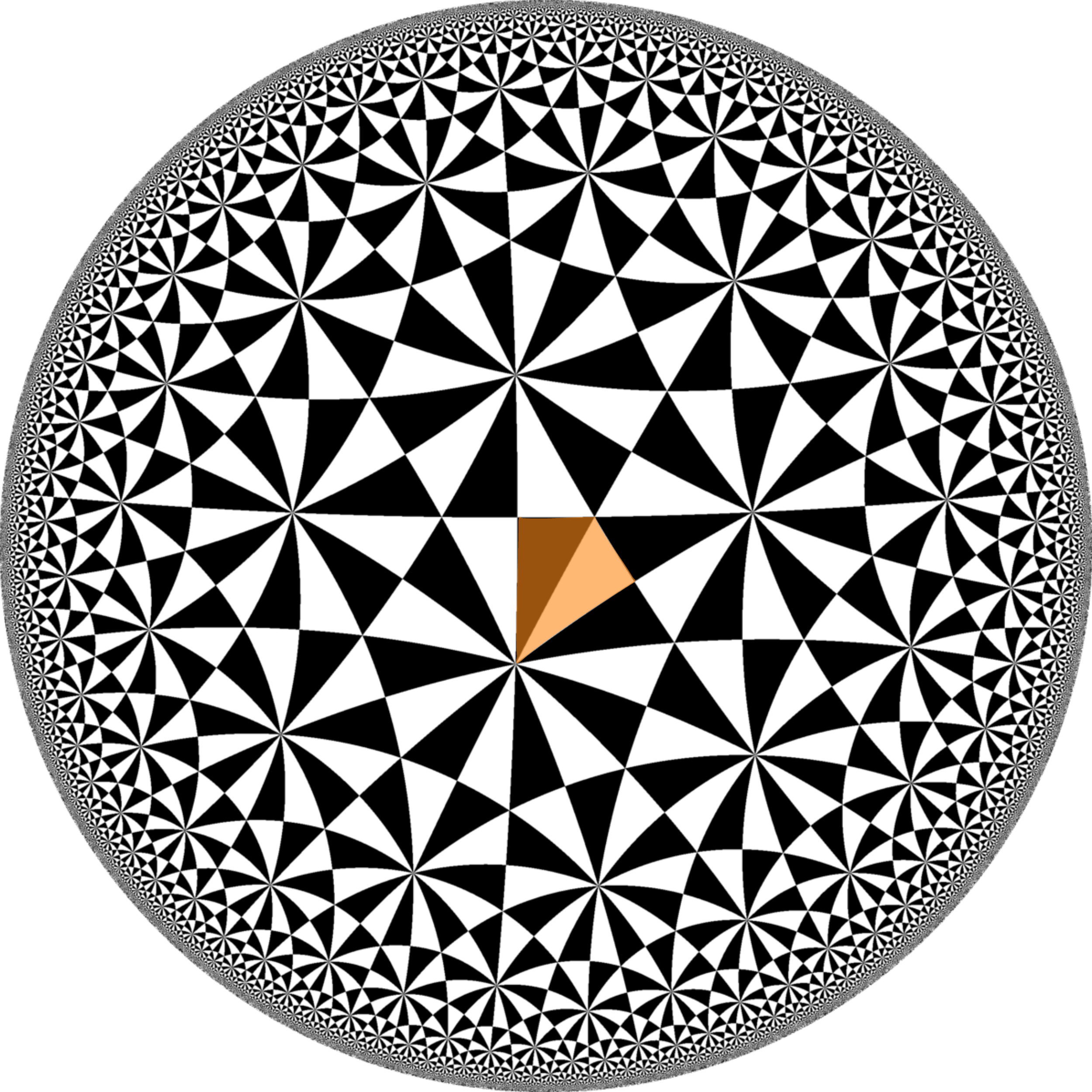}\quad
\includegraphics[width=.3\textwidth]{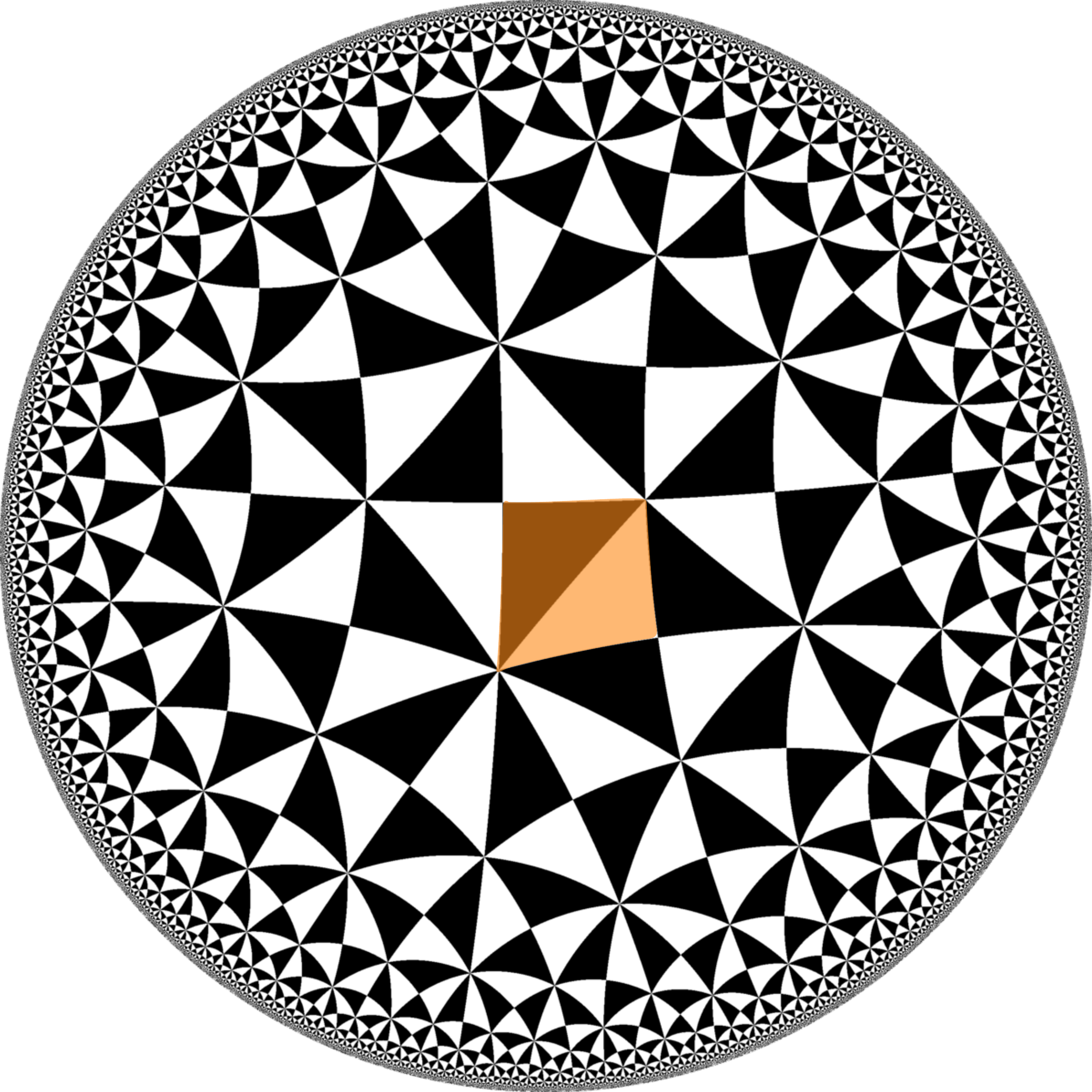}\quad
\includegraphics[width=.3\textwidth]{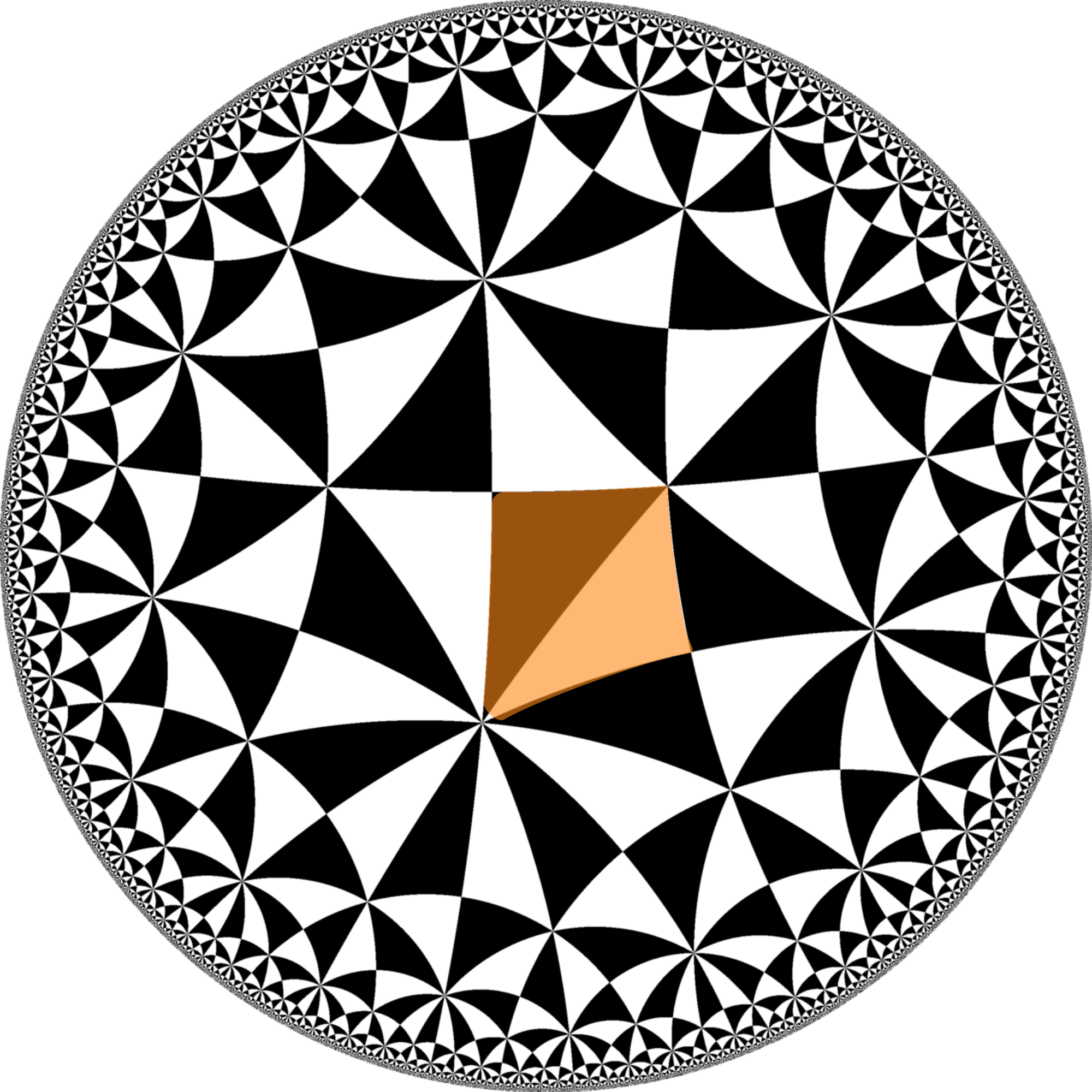}

\includegraphics[width=.3\textwidth]{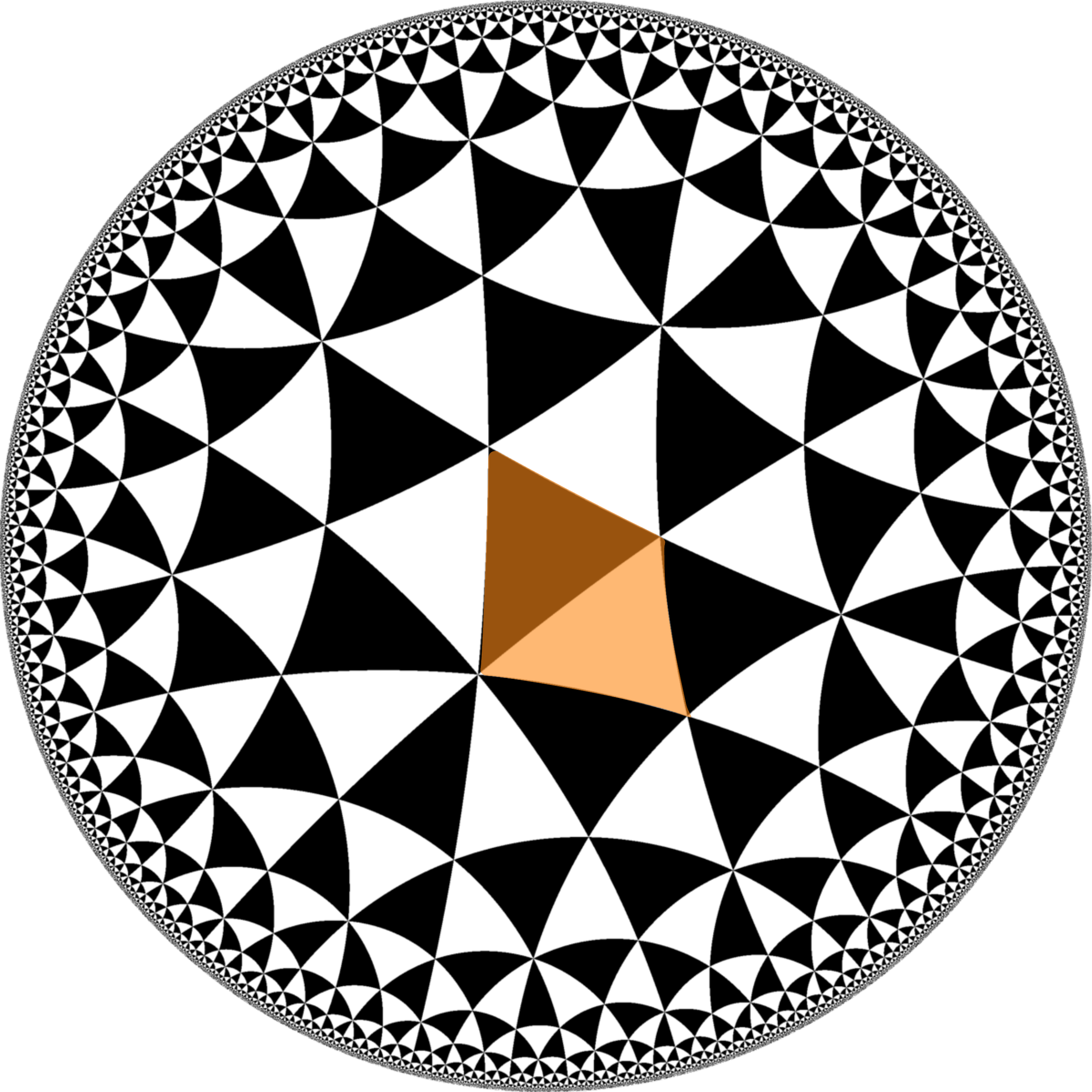}\quad
\includegraphics[width=.3\textwidth]{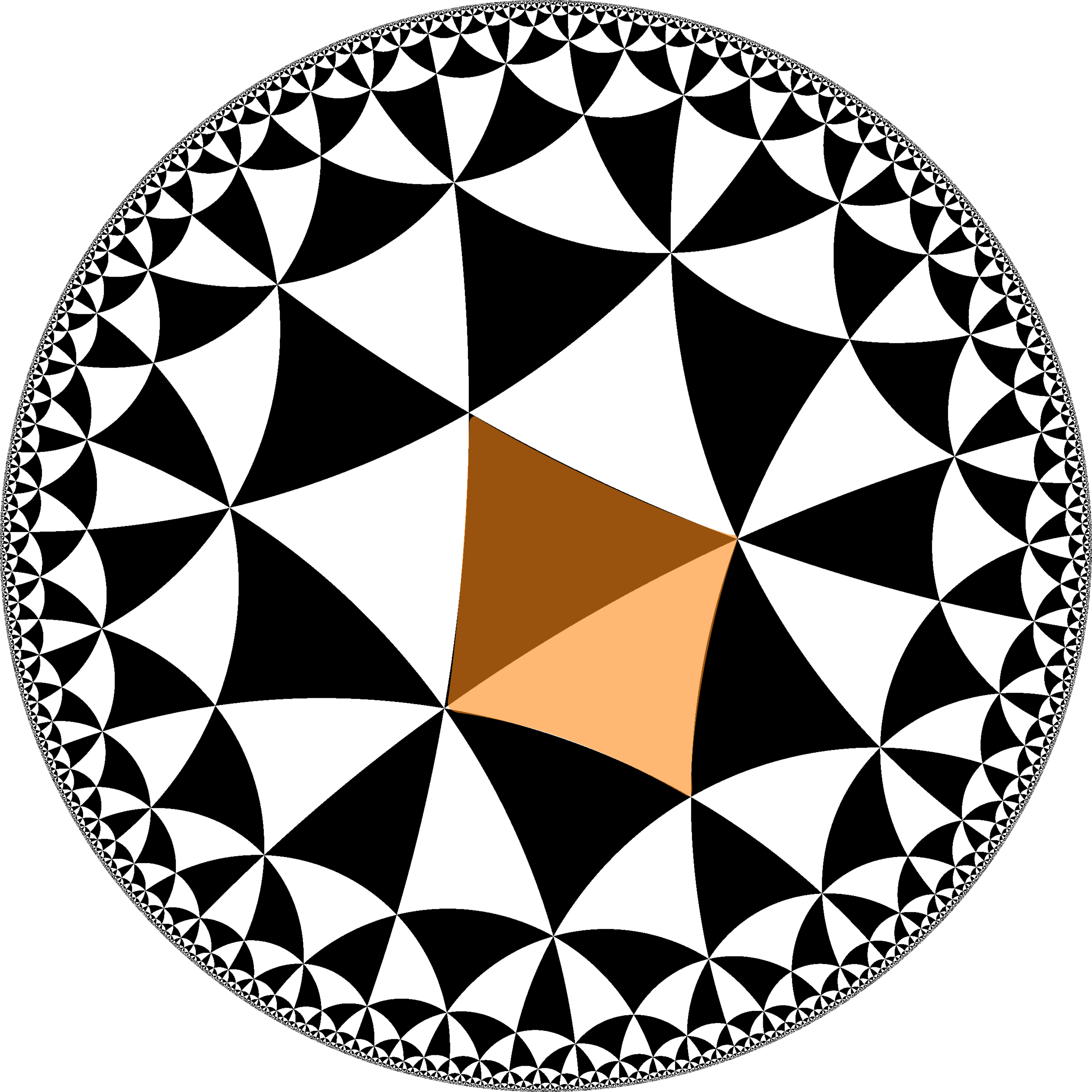}
\end{center}

The quotient~$\Hy/\Gpqr$ is a \term{hyperbolic triangular orbifold}; we denote it by~$\Spqr$. 
It is homeomorphic to a sphere. 
Metrically all points but the projections of~$P, Q, R$ admit neighborhoods isometric to hyperbolic discs. 
The projections of $P, Q, R$ admit neighborhoods isometric to the quotient of a hyperbolic disc by rotations of respective order $p, q, r$. 
%
%

The unit tangent bundle to the hyperbolic plane~$\U\Hy$ is the circle bundle over~$\Hy$ made of length 1 tangent vectors. 
It is homeomorphic to~$\D^2\times \Sph^1$. 
Since isometries also act on tangent vectors, the action of~$\Gpqr$ extends to~$\U\Hy$. 
Note that a rotation~$r$ of order~$k$ on a 2-dimensional disc~$D^2$ acts by screw-motion on the unit tangent bundle~$\U \D^2$. 
Thus the quotient~$\U D^2/\langle r\rangle$ is a 3-manifold where the fiber of the cone point has Seifert invariant $(k,1)$: regular fibers wind $k$ times along the singular fiber, and make one meridional turn around it (right).

\begin{center}
\includegraphics[height=5cm]{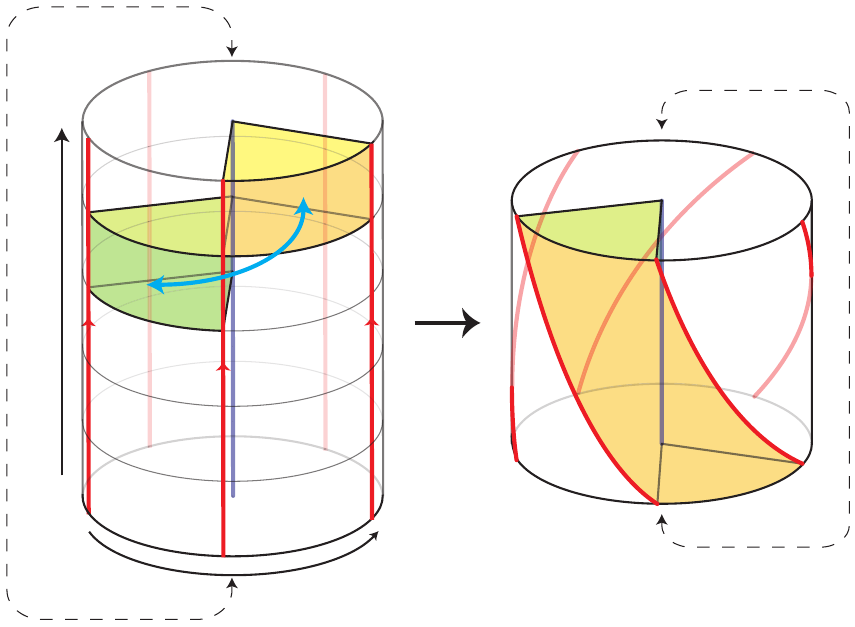}
\end{center}


This implies that the quotient~$\U\Hy/\Gpqr$ is an actual 3-manifold; we denote it by~$\U\Spqr$.
It is Seifert fibered by the fibers of the points of~$\Spqr$.
Regular points have regular fibers, and the three cone points have singular fibers of Seifert invariants~$(p,1), (q,1), (r,1)$ respectively.
We refer to Jose Maria Montesinos' beautiful book~\cite{Montesinos} for more on orbifolds and their unit tangent bundles.
%

\subsection{Geodesic flows and anosovity}

The geodesic flow on~$\U\Hy$ is the flow whose orbits are tangent vectors to directed geodesics, travelled at unit speed. 
Namely, if $\gamma(t)$ denotes a geodesic travelled at speed~$1$, the geodesic flow at time~$t$ sends $(\gamma(0), \dot\gamma(0))$ to $\phi^t((\gamma(0), \dot\gamma(0))=(\gamma(t), \dot\gamma(t))$. 
The action of the geodesic flow commutes with the action of any Fuchsian group, so it descends to a well-defined flow~$(\fpqr)_{t\in\R}$ on~$\U\Spqr$, also called the \term{geodesic flow}. 
Its orbits are lifts of geodesics of~$\Spqr$. 
Note that when a geodesic on~$\Spqr$ goes into a cone point of even order, it makes a U-turn. 
When it goes through a cone point of odd order, it continues straight. 

The geodesic flows on compact quotients of~$\U\Hy$ are the oldest known examplex of Anosov flows~\cite{Anosov}. 
We give here the topological definition which is better suited to our purpose. 
It was only recently shown by Mario Shannon that, for transitive flows in dimension 3, this definition yields the same flows as the classical definition~\cite{Shannon}. 
A 3-dimensional flow is~\term{topologically Anosov} if there are two transverse 2-dimensional foliations $\F^s$ and~$\F^u$ that are invariant under the flow, and such that in every leaf of~$\F^s$ ({\it resp.} $\F^u$) orbits converge exponentially fast in positive ({\it resp.} negative) time. 
The foliations $\F^s$ and $\F^u$ are called \term{stable} and \term{unstable foliations} respectively. 
Note that being topologically Anosov only depends on the underlying 1-dimensional foliation, but not of the actual time-parametrization of the orbits. 

\begin{center}
\includegraphics[width=.5\textwidth]{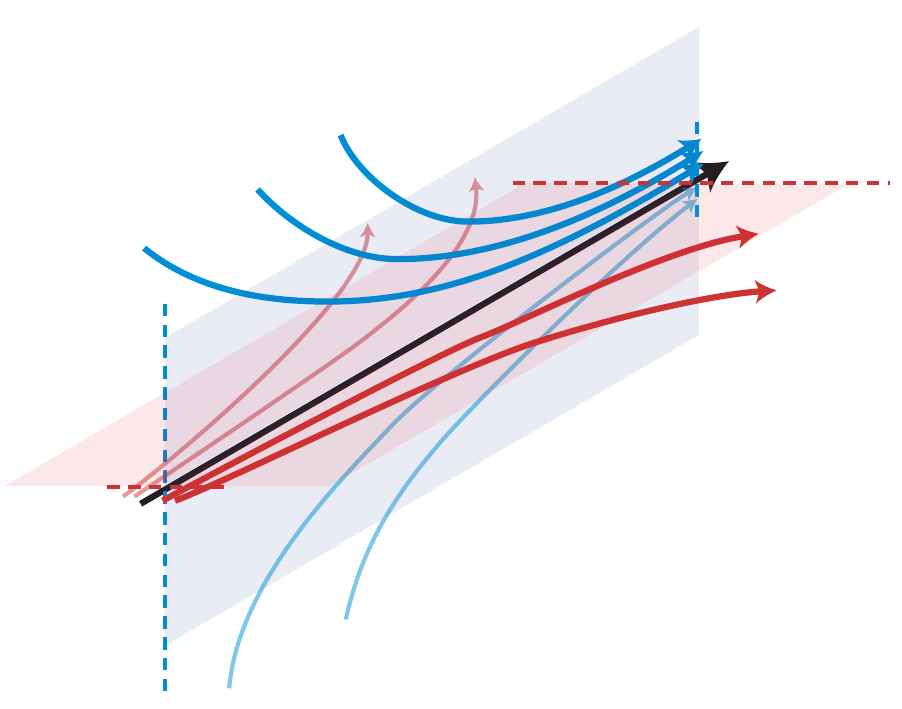}
\end{center}

In the case of the geodesic flow on a quotient of the hyperbolic plane, $\F^s$ corresponds to the vectors tangent to those geodesics that have the same positive extremity in~$\Hy$, while $\F^u$ corresponds to the vectors tangent to those geodesics that have the same negative extremity. 

\begin{center}
\includegraphics[height=4cm]{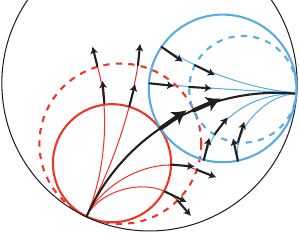}
\end{center}


In general, there is a dichotomy depending on whether the invariant foliations are coorientable or not. 
All the flows we consider here have coorientable invariant foliations. 
This implies in particular that following a periodic orbit, one sees two stable and two unstable half-leaves on each side, and these never get permuted when following the orbit. 

\subsection{Birkhoff sections, first-return maps, and (pseudo-)Anosov maps}
\label{S:Birkhoff}

For $M$ a compact 3-manifold, and $(\phi^t)_{t\in \R}$ a flow with no fixed point on~$M$, a \term{Birkhoff section} for~$(M, \phi^t)$ is the image of an immersed compact surface with boundary~$i:(S, \partial S)\to M$ such that 
\begin{enumerate}
\item the restriction of~$i$ to the interior of~$S$ is an embedding and its image is transverse to the orbits of~$\phi^t$;
\item the restriction of~$i$ to the boundary of~$S$ is a submersion on finitely many periodic orbits;
\item every orbit of~$\phi^t$ intersects~$S$ in bounded time. 
\end{enumerate} 

\begin{center}
\includegraphics[width=.5\textwidth]{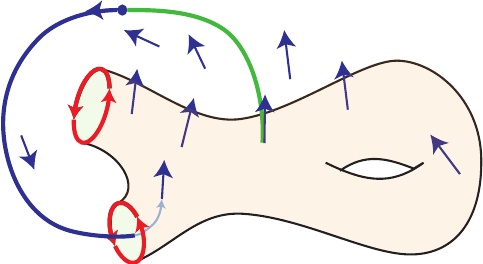}
\end{center}

Note that (2) implies that several boundary components may be mapped on the same periodic orbit of~$\phi^t$. 
Also it allows a given boundary component to wrap several times around a periodic orbit in the longitudinal direction. 

In order to better study the behavior around the boundary, it is convenient to blow-up the boundary orbits as follows: 
given a Birkhoff section $S$, denote by~$M_{\partial S}$ the 3-manifold with toric boundary obtained by blowing up the (finitely many) periodic orbits of~$\phi^t$ on which the boundary of~$S$ maps. 
Each boundary component of~$M_{\partial S}$ is a torus that corresponds to the blow-up of a periodic orbit. 
Topologically, $M_{\partial S}$ is the same as the complement of a tubular neighborhood of~$\partial S$. 
However, assuming~$\phi^t$ to be smooth, it extends to a smooth flow on~$\partial M_{\partial S}$, denoted by~$\hat\phi^t$.
The interior of~$S$ then embeds into~$M_{\partial S}$. 
The Birkhoff section~$S$ is called~\term{$\partial$-strong}\footnote{This refinement was proposed by Umberto Hryniewicz to formalize an assumption that was often done but never stated.} if $S^\circ$ extends into a compact surface with boundary~$\hat S$ in~$M_{\partial S}$ such that every boundary component of~$\hat S$ is a curve in~$\partial M_{\partial S}$ that is transverse to the extended flow~$\hat\phi^t$. 
Note that a Birkhoff section can always be perturbed into a $\partial$-strong Birkhoff section. 

Starting from a point in~$\hat S$ and following~$\hat\phi^t$ there is a first-return time on~$\hat S$ and a induced \term{first-return map}~$f_{\hat S}:\hat S\to \hat S$ which is a homeomorphism. 
Blowing down every boundary component of~$\hat S$ into a point, we obtain a compact surface~$\bar S$ with no boundary, and a map~$f_{\bar S}: \bar S\to\bar S$, which is still a homeomorphism\footnote{In order to describe~$\bar S$ and the first-return map~$f_{\bar S}$, one does not really need $\partial$-strongness and looking at~$\hat S$, but we find it convenient to explain it in this setting.}. 

Let $\gamma$ be a periodic orbit of~$\phi^t$ on which lie $c$ boundary components of~$S$. 
Then the boundary torus of~$M_{\partial S}$ corresponding to $\gamma$ contains~$c$ components of~$\partial \hat S$ which are parallel. 
These $c$ curves then give $c$ points in~$\bar S$, which form a periodic orbit for~$f_{\bar S}$ of period~$c$. 

\bigskip


When $\phi^t$ is Anosov with coorientable invariant foliations, there are two canonical directions on the boundary components of~$\partial M_{\partial S}$ given by the trace of the stable/unstable directions and the meridians. 
These directions allow to describe every boundary curve of~$\hat S$ by a direction $(a,b)$ of coprime integers with $b>0$, where $a$ is the (algebraic) intersection with any meridian curve and $b$ is the intersection with the trace of one stable half-leaf (on the following picture one has~$(a,b)=(2,1)$, so that the slope $b/a$ is $\frac12$).

\begin{center}
\includegraphics[width=.5\textwidth]{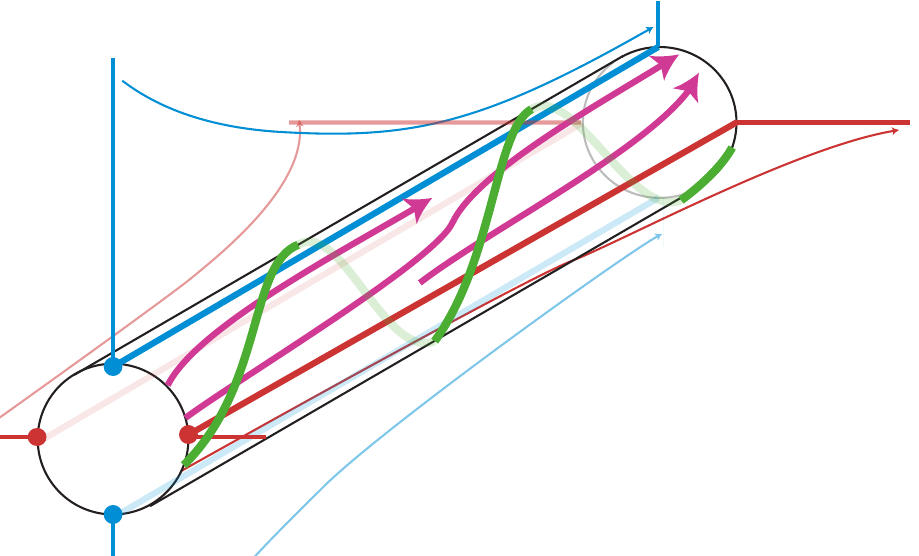}
\end{center}

In particular, assuming $S$ to be $\partial$-strong, every boundary component of~$\hat S$ is a curve on~$\partial M_{\partial S}$ whose direction is of the form $(a,b)$ with $b> 0$. 
Indeed, $b=0$ would violate condition~(3): the boundary curve would be parallel to the stable and unstable directions, and the first-return time would not be bounded. 

When~$\phi^t$ is Anosov, one can intersect the invariant foliations of~$\phi^t$ with a $\partial$-strong section~$S$. 
This yields two 1-dimensional foliations on $S^\circ$, that are invariant by the first-return map. 
These foliations actually extend to~$\hat S$ into foliations with boundary with the local form depicted below.

\begin{center}
\includegraphics[width=.5\textwidth]{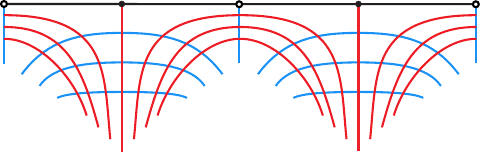}
\end{center}

We denote these by~$\F_{\hat S}^s$ and $\F_{\hat S}^u$. 
One checks that if $s$ is a boundary component of~$\hat S$ with direction~$(a,b)$, then the extended foliations have $2b$ singular half-leaves on~$s$. 


The case that is of interest in this paper is when $b=1$ for every boundary component of~$\hat S$. 
In this case only, the foliations~$\F_{\hat S}^s$ and $\F_{\hat S}^u$ project to 1-dimensional foliations without singularities on~$\bar S$. 
This implies first that $\bar S$ is a (orientable) compact surface with empty boundary that carries two non-singular foliations, hence it is a torus. 
Second the first-return map~$f_{\bar S}$ preserves two non-singular foliations, expanding the first one and contracting the second, so it is an Anosov map of the torus. 
It is well-known that such a map is conjugated to the action of a linear map~$A\in\SLZ$ with $\tr(A)>2$, and that two such maps~$A$ and $B$ are conjugated if and only if $A$ and $B$ are conjugated in~$\SLZ$~\cite[Thm 18.6.1]{KH}.

Denoting by~$X$ the matrix~$(\begin{smallmatrix} 1&1\\0&1\end{smallmatrix})$ and by $Y$ the matrix~$(\begin{smallmatrix} 1&0\\1&1\end{smallmatrix})$, conjugacy classes in~$\SLZ$ correspond to words in~$X$ and $Y$ containing both letters, up to cyclic permutations. 
In particular, there is only one conjugacy class of such maps with one fixed point on~$\R^2/\Z^2$ only, namely the conjugacy class of the matrix~$XY=(\begin{smallmatrix} 2&1\\1&1\end{smallmatrix})$.

Summarizing the above discussion one has

\begin{proposition}
Assume that $M$ is a 3-manifold that supports an Anosov flow~$\phi^t$, that $\gamma$ is a periodic orbit of~$\phi^t$ which is the boundary of a Birkhoff section~$S$ so that the boundary direction is~$(a,b)$. 
Then the manifold~$M(\gamma, b/a)$ is the mapping torus of the first-return map on~$\bar S$ along~$\phi^t$. 
Moreover, when $b= 1$, the surface~$\bar S$ is necessarily a torus.
\end{proposition}


\section{The $237$-case}

Here $p=2, q=3, r=7$, so that the triangle~$PQR$ has angles $\pi/2, \pi/3, \pi/7$. 
In the 2-dimensional orbifold~$\S_{237}$ we consider the geodesic~$h$ that is obtained by lifting in~$\S_{237}$ the altitude of~$P$ in~$PQR$. 
Denote by~$\vec h$ its lift in~$\U\S_{237}$.

\begin{center}
\includegraphics[height=4.5cm]{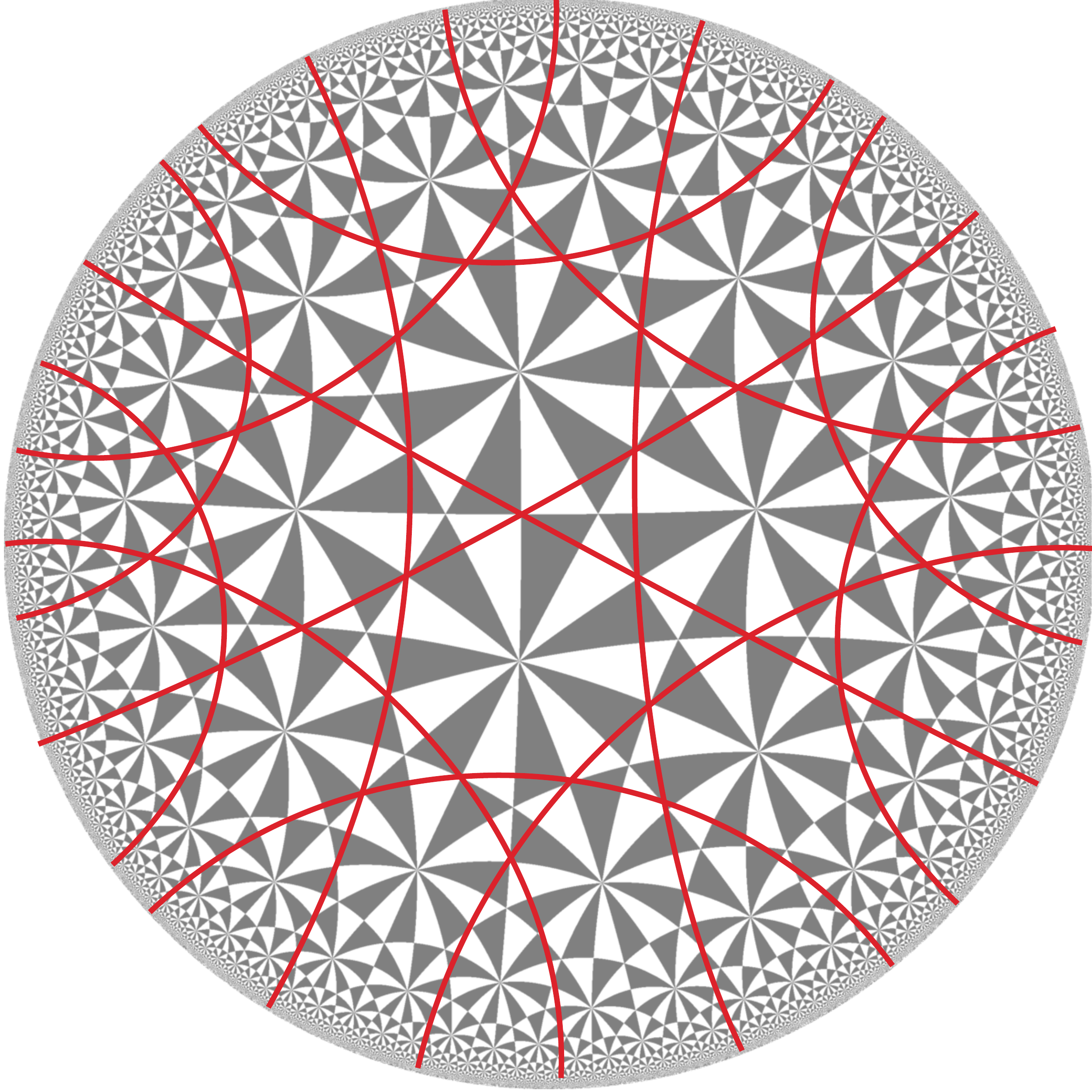}
\qquad\qquad
\includegraphics[height=4.5cm]{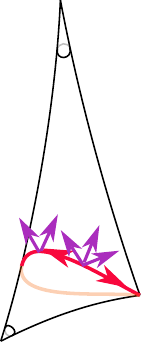}
\qquad\qquad
\includegraphics[height=4cm]{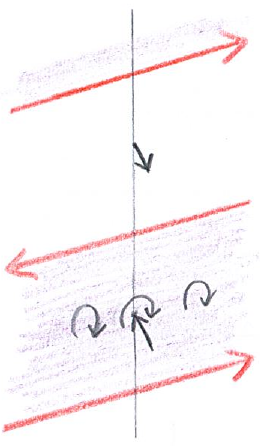}
\end{center}

Since $h$ goes through the point~$P$ which has order~2, there is one lift only in~$\U\S_{237}$, which consists in following $h$ with both orientations.
Denote by~$S_{237}$ the 2-chain in~$\U\S_{237}$ which consists in all tangent vectors based on~$h$ and pointing to the side of~$R$. 
It is cooriented by the orbits of the geodesic flow~$\phi^t_{237}$. 
With the induced orientation, its oriented boundary is~$-\vec h$. 

%

The next picture shows how $S_{237}$ looks around the fiber of the order 2-cone point. 
Combinatorially, it is a rectangle with some side identifications. 

\begin{picture}(180,240)(0,0)
\put(60,150){\includegraphics[height=3cm, angle=1]{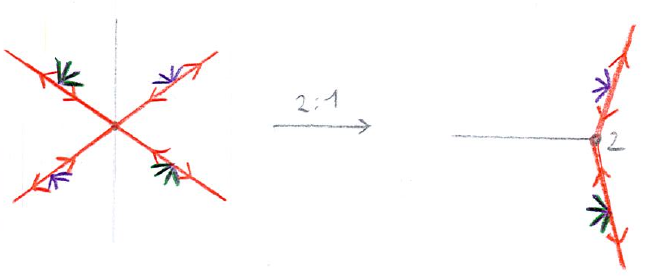}}
\put(25,-20){\includegraphics[height=6.5cm]{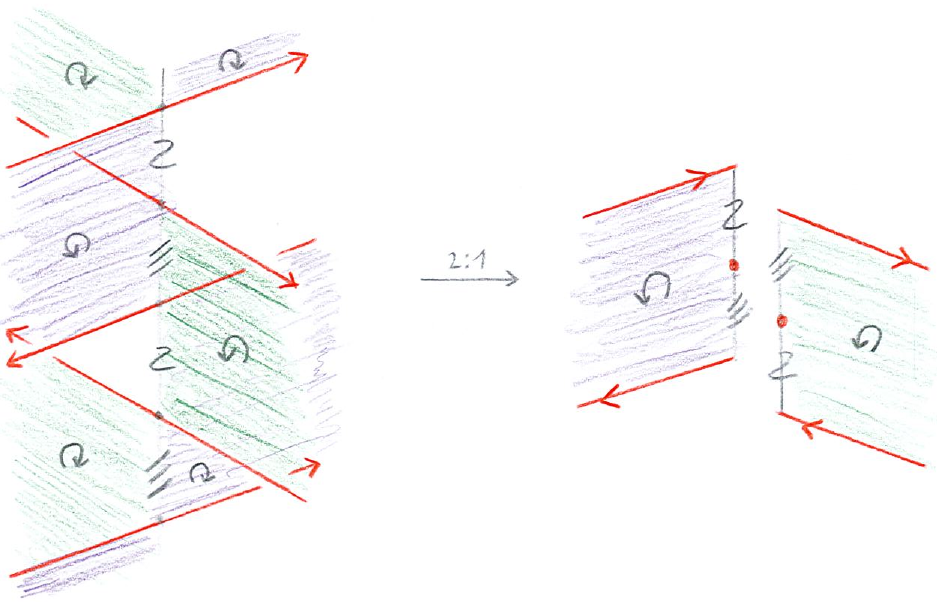}}
\end{picture}

One checks that the Euler characteristic of~$S_{237}$ is~$-1$. 
Thus it is a torus with one boundary component. 
As explained in the end of Section~\ref{S:Prelim}, one can blow down the boundary component of~$S_{237}$ into a point, yielding a surface~$\bar S_{237}$ and a first-return map~$f_{\bar S_{237}}:\bar S_{237}\to\bar S_{237}$ that is Anosov. 

We are left with the computation of the first-return map $f_{\bar S_{237}}$ on~$\bar S_{237}$ along~$\phi_{237}^t$. 
One option is to make an explicit computation, as was done in a more general setting~\cite{D:ETDS, DLiechti}. 
Another (easier?) way is to show that $f_{\bar S_{237}}$ has one fixed point only. 
Here the boundary projects to a fixed point after collapsing. 
One can then look at all other geodesics in~$\Hy$, and convince oneself that no periodic geodesic intersects~$S_{237}$ only once. 
Indeed a geodesic intersecting~$S_{237}$ corresponds to an element~$g\in\Gamma_{237}$ that translates an heptagon bounded by the lifts of~$h$ to an adjacent heptagon, for example the green one to the pink one in the next picture. 

\begin{center}
\includegraphics[height=7cm]{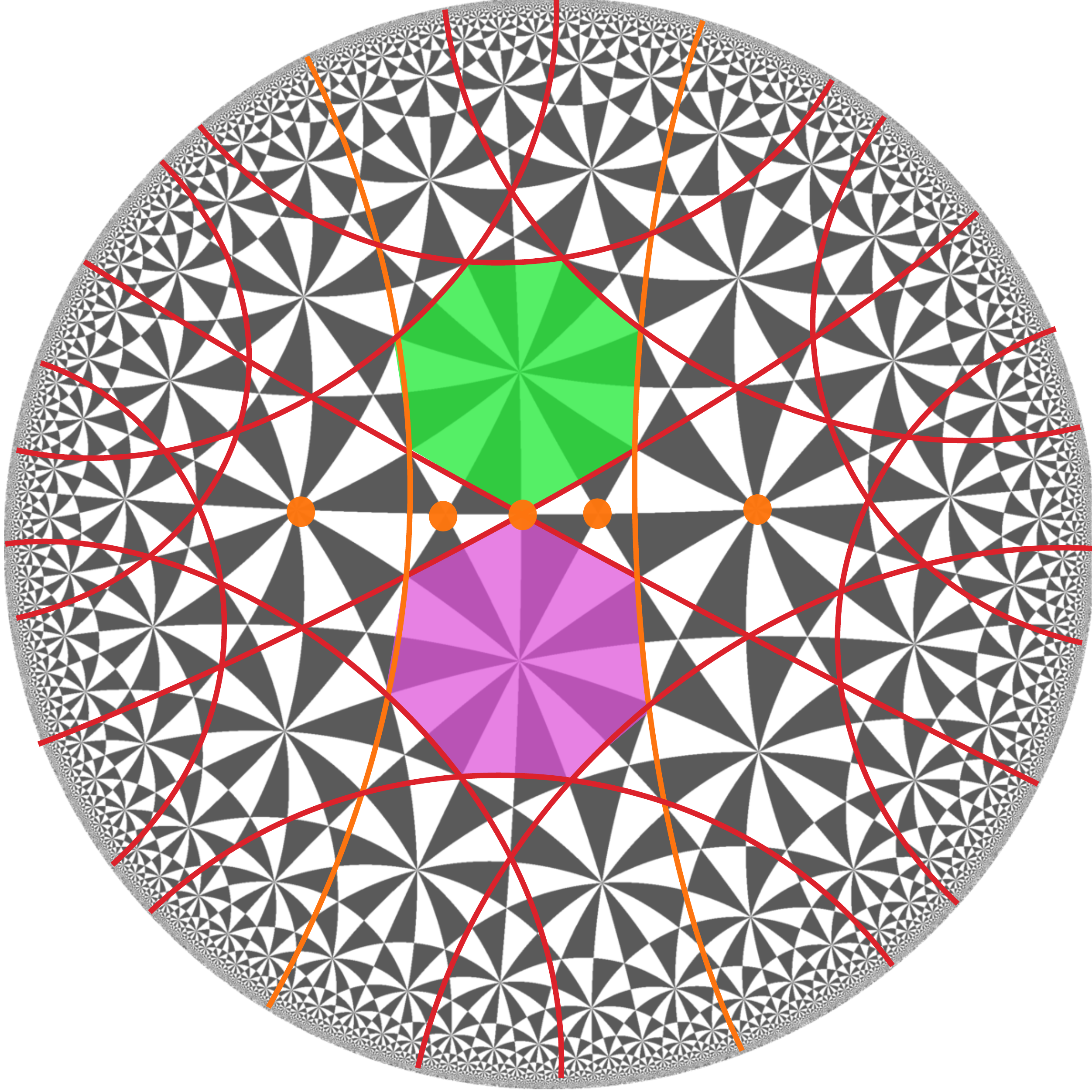}
\end{center}

There are 7 such isometries, 5 of which are rotations (around the orange points), and 2 of which are translation along copies of~$h$. 
Thus there is no fixed point other than the one coming from~$h$. 


\section{The $245$-case}

In the orbifold~$\S_{245}$ we consider the geodesic $b$ that is the shortest segment connecting $P$ and $Q$. 
Denote by $\vec b$ its lift in~$\U\S_{245}$. 
Since $b$ goes though two cone points of even orders, it makes a U-turn at both points, hence there is only one lift for~$b$ in~$\U\S_{245}$. 

\begin{center}
\includegraphics[height=5cm]{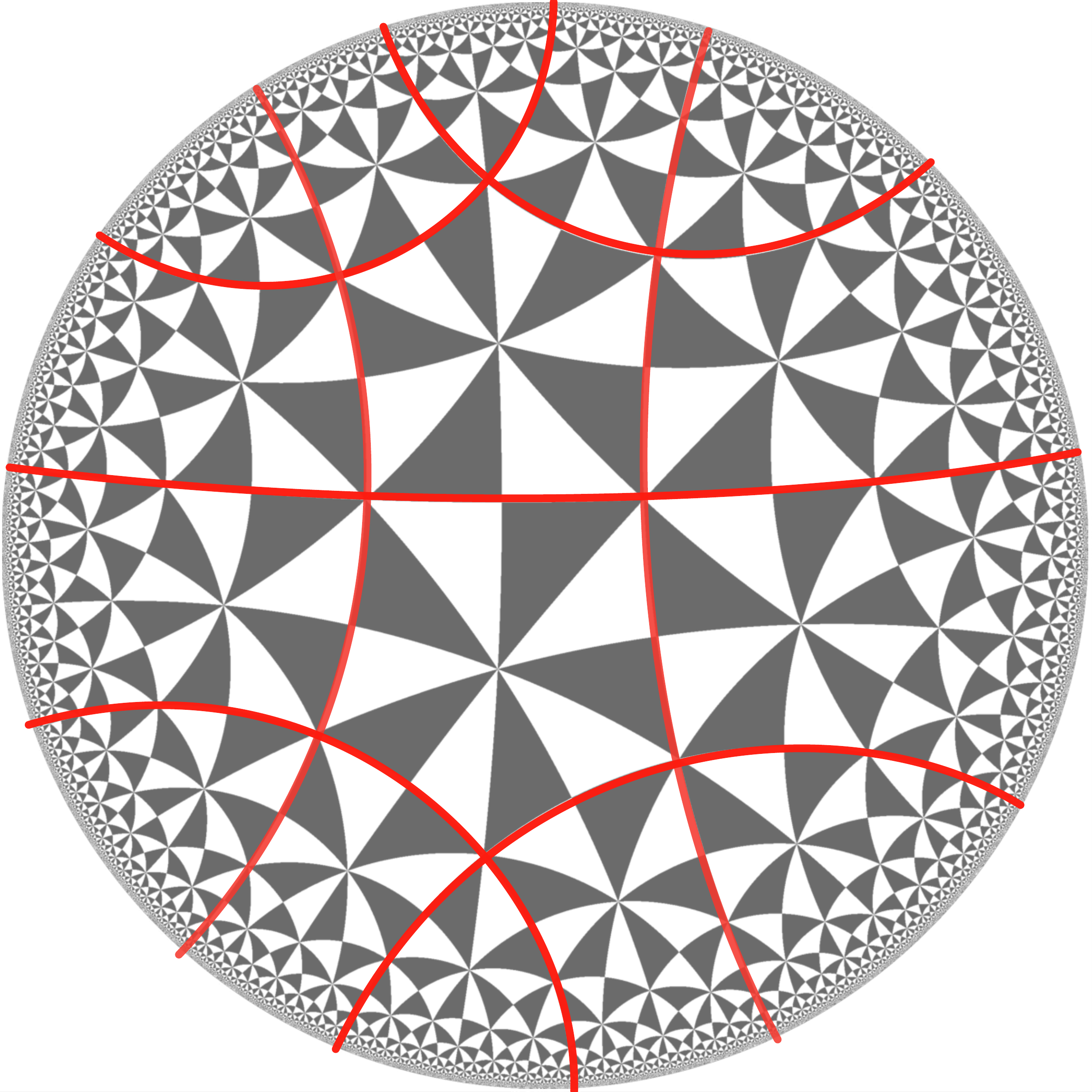}
\qquad
\includegraphics[width=.43\textwidth]{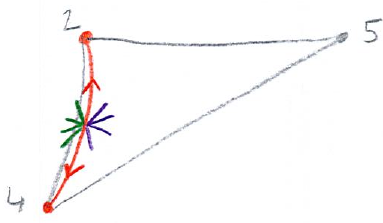}
\end{center}

Denote by~$\hat S_{245}$ the surface in~$\U\S_{245}$ that consists in all tangent vectors based on~$b$. 
It is cooriented by the orbits of~$\phi_{245}^t$. 
With the induced orientation, its oriented boundary is~$-2\vec b$. 

Around the fiber of~$P$, the situation is a bit different from the~$237$-case since $b$ comes at~$P$ from one direction only, but $S_{245}$ consists of two ribbons.  
The local picture at~$\U P$ is the following. 
The surface~$\hat S_{245}$ is not singular here.

\begin{center}
\includegraphics[height=3.5cm]{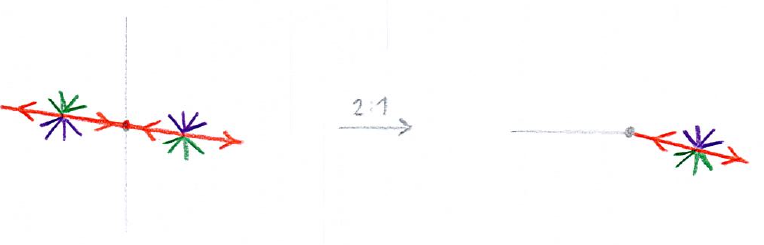}

\includegraphics[height=5.5cm]{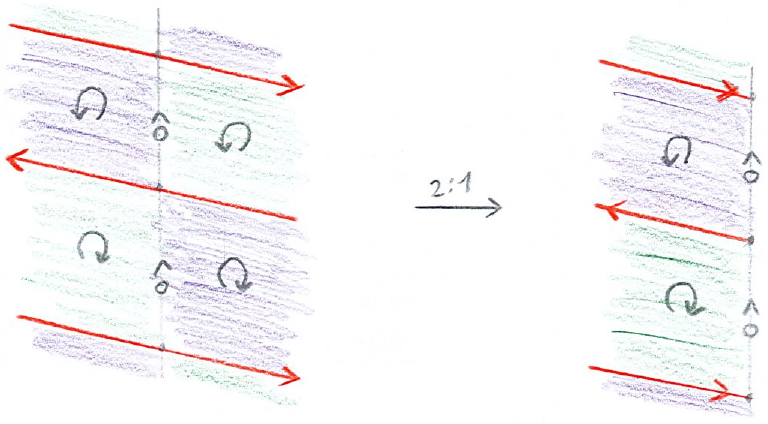}
\end{center}

The situation around the fiber of~$Q$ is more involved and the surface~$\hat S_{245}$ is singular there. 
In a local chart that is a degree 4-cover, the lift of~$\hat S_{245}$ consists of 8 pieces: 2 above each segment that covers~$b$. 
Coloring in green those pieces that correspond to vector oriented clockwisely around~$Q$ and coloring in purple the other four pieces, the order-4 rotation identifies the 4 green pieces together and the 4 purple pieces together. 

\begin{center}
\includegraphics[height=3cm]{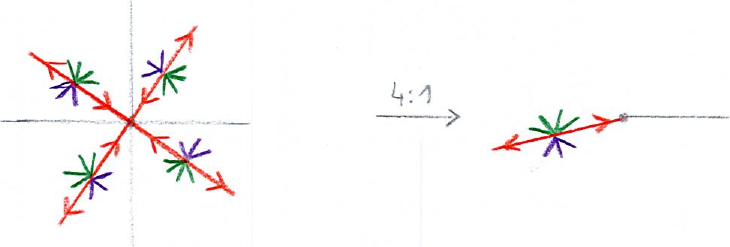}

\includegraphics[height=7cm]{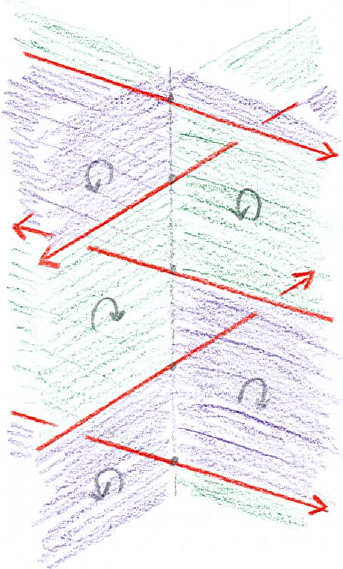}
\hspace{1cm}
\includegraphics[height=7cm]{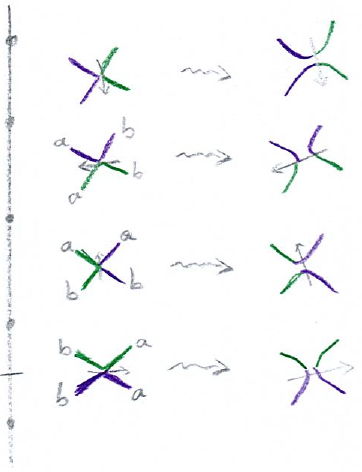}

\includegraphics[height=5cm]{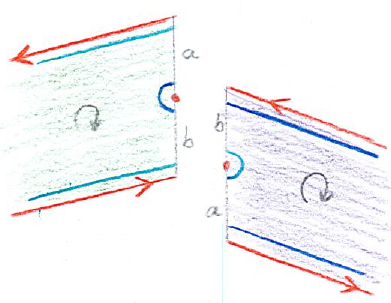}
\end{center}

One can desingularize $\hat S_{245}$ transversally to~$\phi_{245}^t$ as shown on the right and obtain a topological surface~$S_{245}$. 
The desingularization connects the top part of the vertical boundary of each green piece with the bottom part of the vertical boundary of another green piece, and similarly with the purple pieces. 

Thus $S_{245}$ is a rectangle with some identifications on the vertical border. 
One checks that its Euler characteristic is~$-1$. 
Also it has one boundary component: indeed, when travelling on~$S_{245}^\circ$ along its border, one checks that there is one component only: in the fiber of $Q$ the borders of the green and purple part do not mix, but in the fiber of~$P$ they get exchanged. 
Thus $S_{245}$ is a torus with one boundary component. 
One can compute the boundary direction of~$S_{245}$ along its unique boundary component: its winds twice along~$\vec b$, and once meridionally, so the direction is~$(2,1)$. 

In order to understand the induced first-return map~$f_{\bar S_{245}}$, we show that is has one fixed point only. 
Such a fixed point would correspond to an element of~$\Gamma_{245}$ mapping a pentagon in~$\Hy$ delimited by the lifts of~$b$ to an adjacent pentagon. 

\begin{center}
\includegraphics[height=7cm]{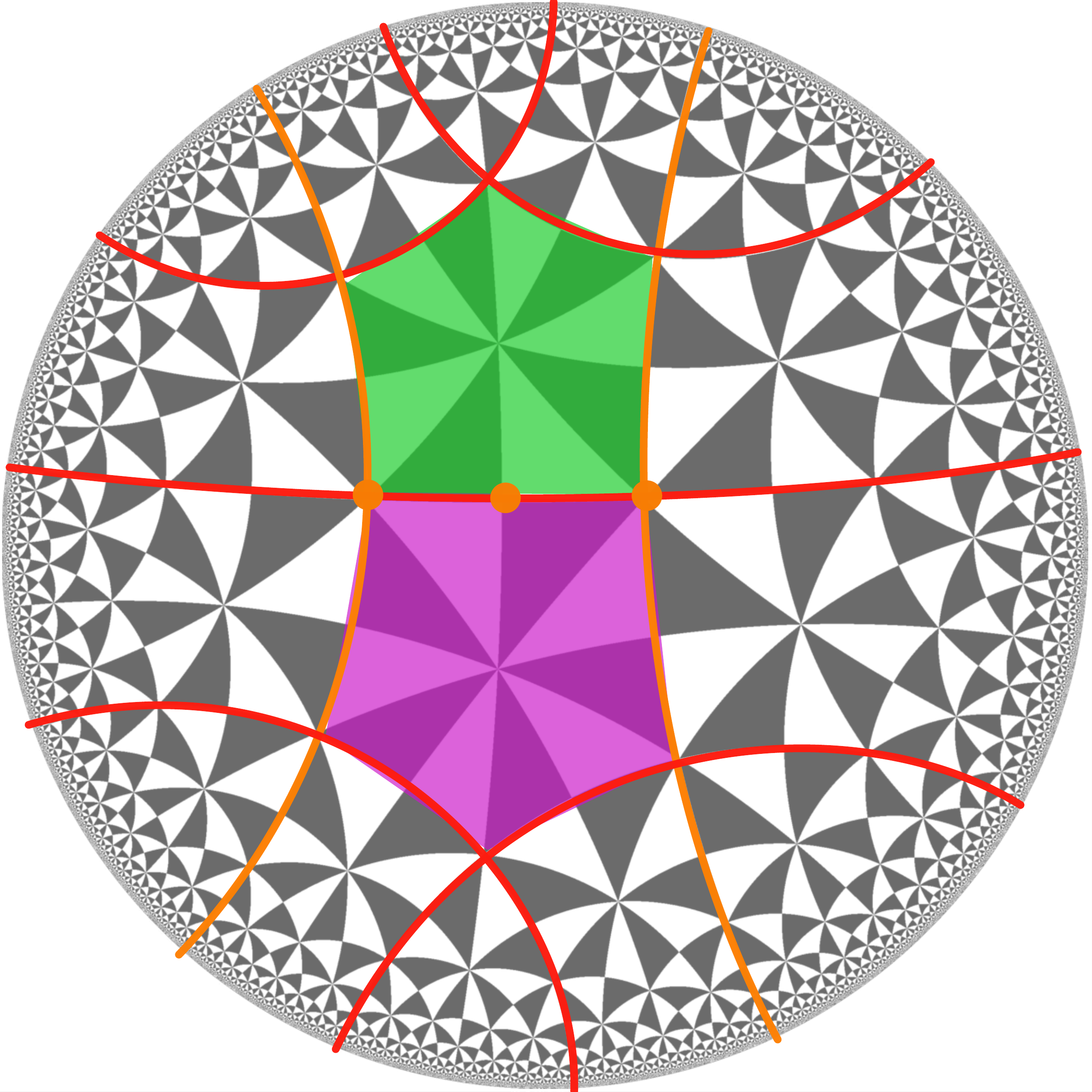}
\end{center}

There are only 5 ways to do so, 3 of which are rotations (around the orange points), and 2 of which are translations along copies of~$b$ (in orange too). 
Thus the first-return map on~$S_{245}$ is also conjugated to the cat-bat map. 

\newpage
\section{The $246$-case}

This case is close to the previous one. 
In the orbifold~$\S_{246}$ we consider the geodesic $c$ that is the shortest segment connecting $P$ and $R$. 
Denote by $\vec c$ its lift in~$\U\S_{246}$. 
Since $c$ goes though two cone points of even points, there is only on lift for~$c$. 

\begin{center}
\includegraphics[height=5cm]{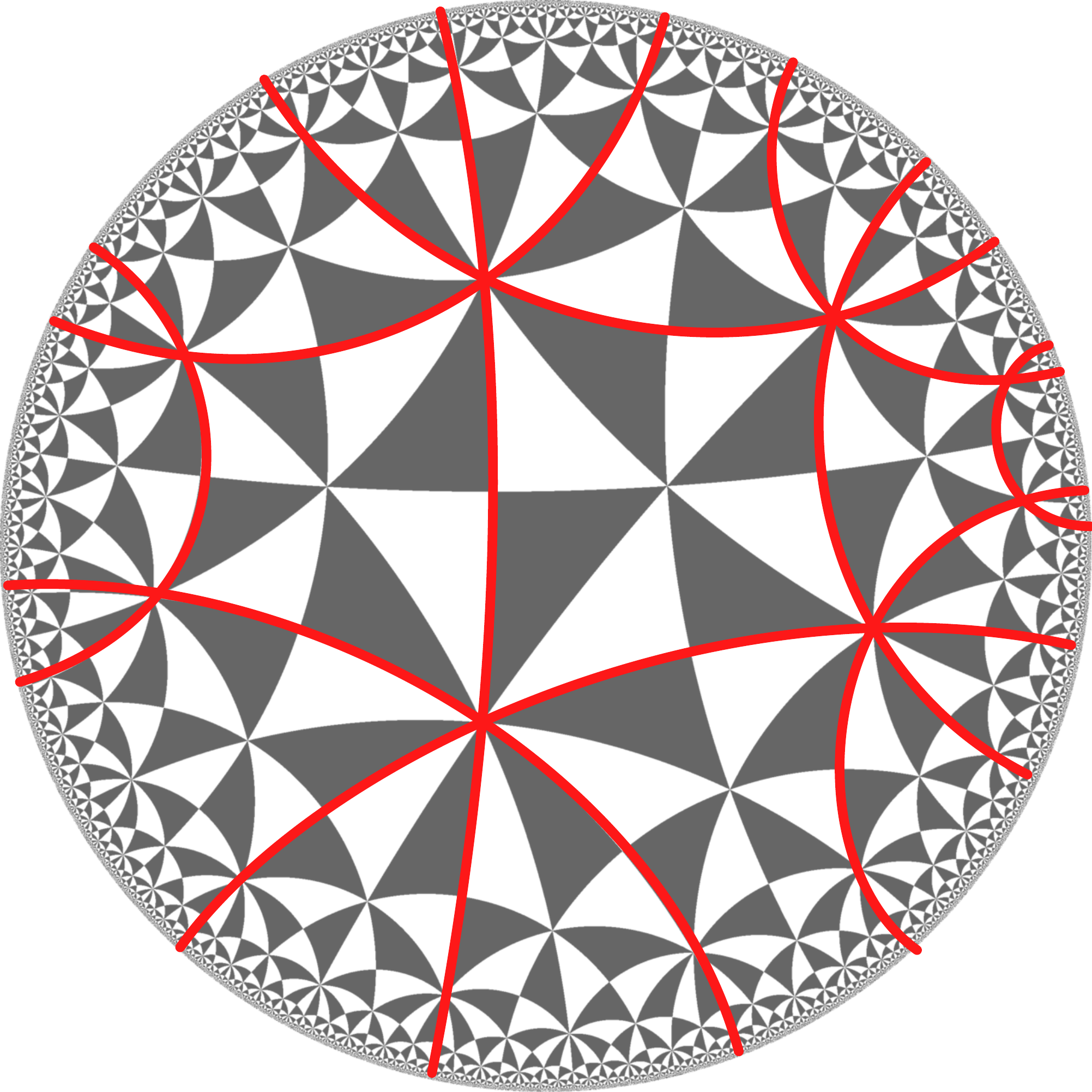}
\qquad
\includegraphics[width=.35\textwidth]{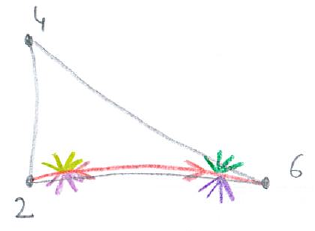}
\end{center}

Denote by~$\hat S_{246}$ the surface in~$\U\S_{246}$ which consists in all tangent vectors based on~$c$. 
It is cooriented by the orbits of~$\phi_{246}$. 
With the induced orientation, its oriented boundary is~$-2\vec c$. 
As in the previous case, $\hat S_{246}$ is not singular around the order 2-cone point. 

\begin{center}
\includegraphics[height=4cm]{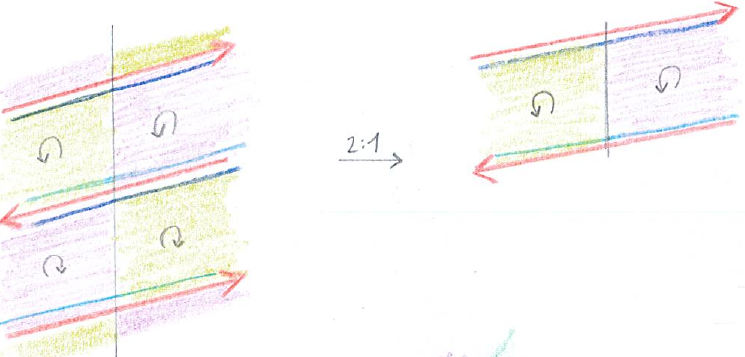}
\end{center}

However the situation around the fiber or~$R$ is more complicated than in the previous case. 
In a local chart that is a degree 6-cover, the lift of~$\hat S_{246}$ consists of 12 pieces: 2 above each segment that covers~$c$. 
Coloring in green those pieces that correspond to vectors oriented clockwisely around~$R$ and coloring in purple the other six pieces, the order-6 rotation identifies all 6 green pieces together and all 6 purple pieces together. 

\begin{center}
\includegraphics[height=7cm]{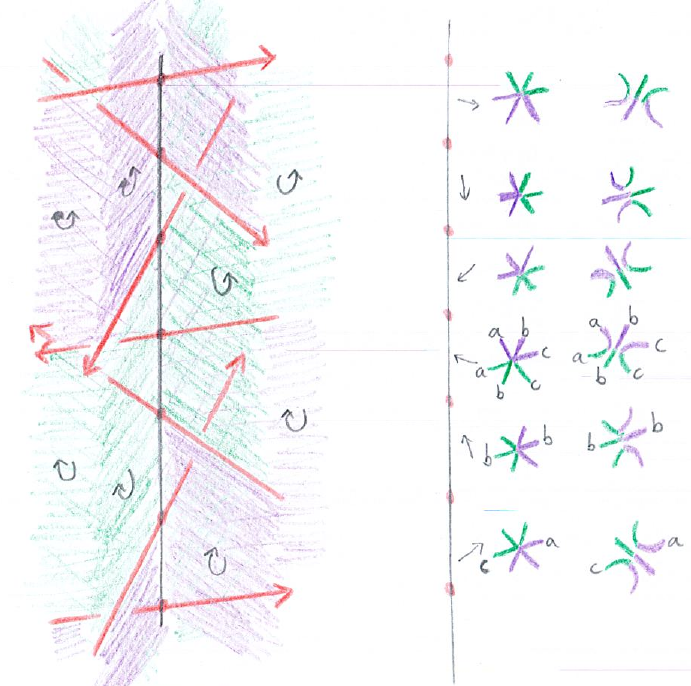}

\includegraphics[height=4cm]{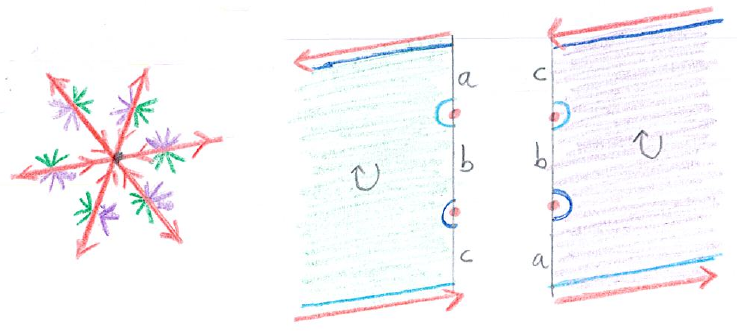}
\end{center}

One can desingularize $\hat S_{246}$ transversally to~$\phi_{246}^t$ as shown on the picture (top right) and obtain a topological surface $S_{246}$. 
The desingularization connects the top part of the vertical boundary of each purple piece with the bottom part of the vertical boundary of an adjacent green piece, the middle part of the vertical boundary of each purple piece with the middle part of the vertical boundary of the opposite green piece, and the bottom part of the vertical boundary of each purple piece with the top part of the vertical boundary of an adjacent green piece (bottom right). 

Thus $S_{246}$ is a rectangle with some identifications on the vertical border. 
One checks that its Euler characteristic is~$-2$. 
Also, travelling along~$\vec c$, one sees that there are two boundary components (light blue and dark blue on the bottom right picture). 
Thus $S_{246}$ is a torus with two boundary components, each with direction~$(1,1)$. 

In order to understand the induced first-return map, we show that is has one fixed point only. 
Unlike the two previous cases, the boundary is not a fixed point: it is a period 2-point since the two boundary components are exchanged by the first-return map. 

\begin{center}
\includegraphics[height=7cm]{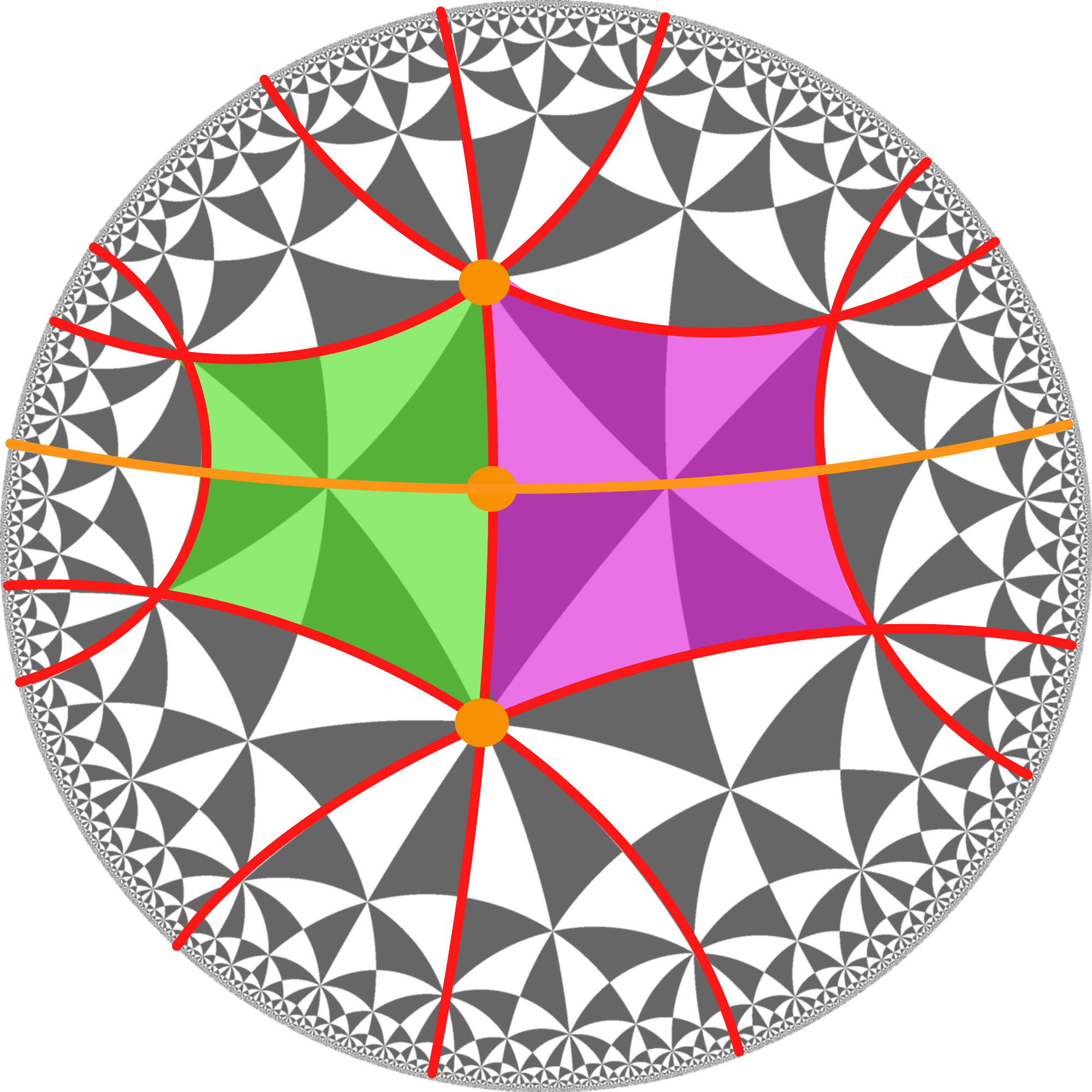}
\end{center}

On the other hand, there is one fixed point in the interior of~$S_{246}$. 
Indeed, such a fixed point corresponds to an element $g\in\Gamma_{246}$ transporting one square in~$\Hy$ bounded by copies of~$c$ (the green square, say) on an adjacent one (the pink one, say). 
One sees that 4 isometries transport the green square on the pink one, namely 3 rotations and 1 translation along the orange geodesic (which actually corresponds to the line $PQ)$. 
The latter yields the unique fixed point. 
Thus the first-return map on~$S_{246}$ is also conjugated to the cat-bat map. 

\section{The $334$-case}

This case and the next one are a bit different since the surfaces we consider are not vertical, but rather horizontal, that is, transverse almost everywhere to the fibering by circles of the unit tangent bundle. 
In that sense, they are closer to the original constructions of surfaces of sections by Birkhoff and~Fried. 
The existence and the properties of the surface~$S_{334}$ we construct here were already proven~\cite{D:TSG}, but the surface was not explicitly described. 
We follow another approach here, reminicent of the surface~$S_{pqr}$ in~\cite{D:ETDS}.

In~$\Hy$ tiled copies of a triangle~$PQR$ with angles $\pi/3, \pi/3, \pi/4$, call $J$ the midpoint of the segment~$PQ$. 
On the orbifold~$\S_{334}$, consider the oriented geodesic~$\gamma_8$ that has one self-intersection at~$J$ and that winds once positively around~$P$ and once negatively around~$Q$. 
In~$\Hy$, the curve~$\gamma_8$ and its images by~$\Gamma_{334}$ tile~$\Hy$ into octagons (containing the order 4-cone points), and two types of triangles containing the orbits under the action~$\Gamma_{334}$ of $P$ and $Q$ respectively. 
Lift~$\gamma_8$ in~$\U\S_{334}$ into an periodic orbit~$\vec\gamma_8$ of the the geodesic flow~$\phi_{334}$. 

\begin{center}
\includegraphics[height=6cm]{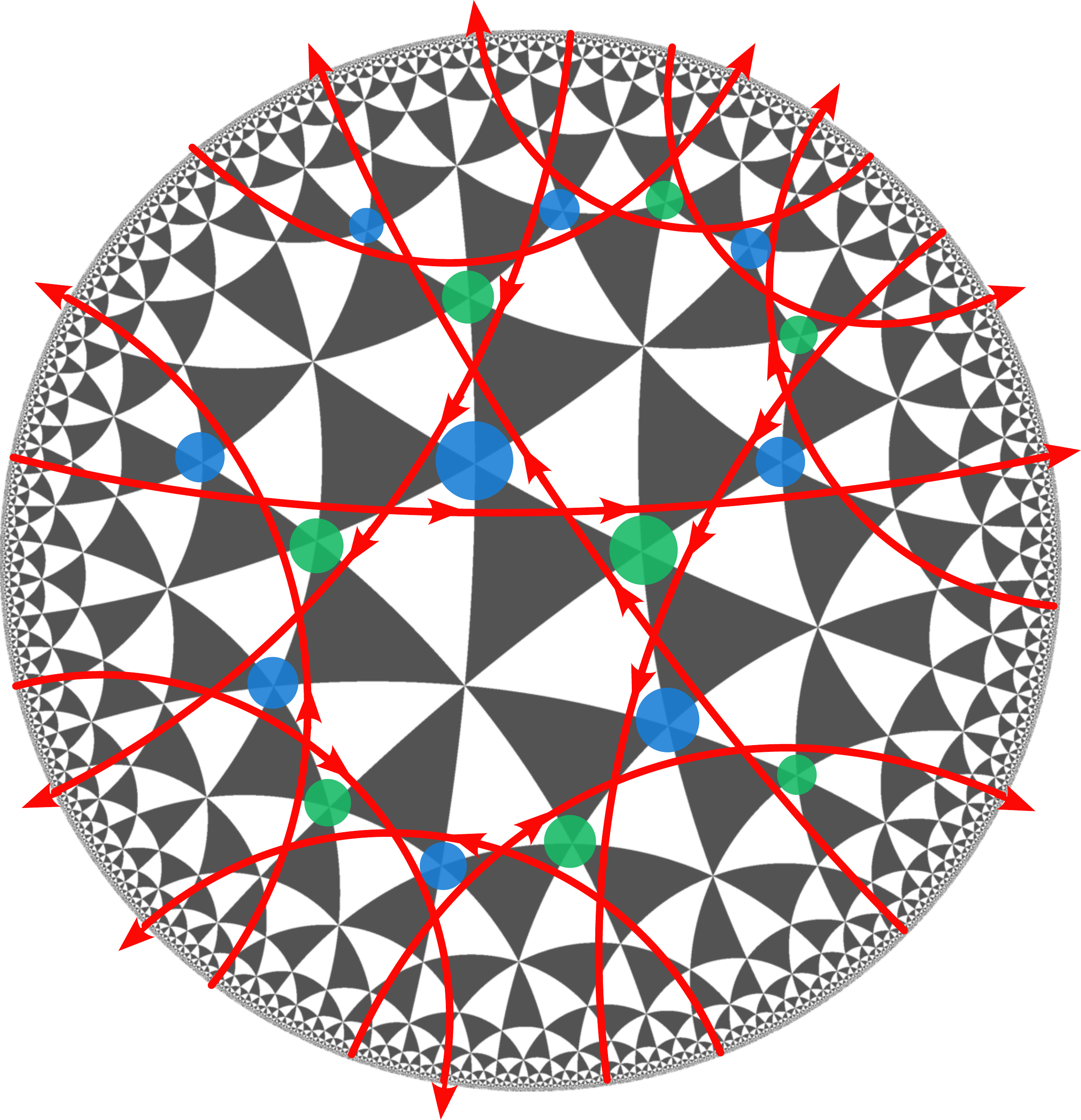}
\qquad
\includegraphics[height=4cm]{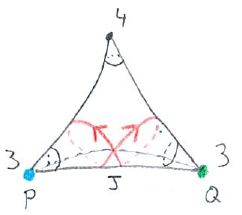}
\end{center}

In~$\Hy$ consider an oriented foliation~$\F_{334}$ of the interiors of the two triangles delimited by~$\gamma_8$ with one vertex at~$J$ by convex curves tangent to~$\gamma_8$ at~$J$. 

\begin{center}
\includegraphics[height=7cm]{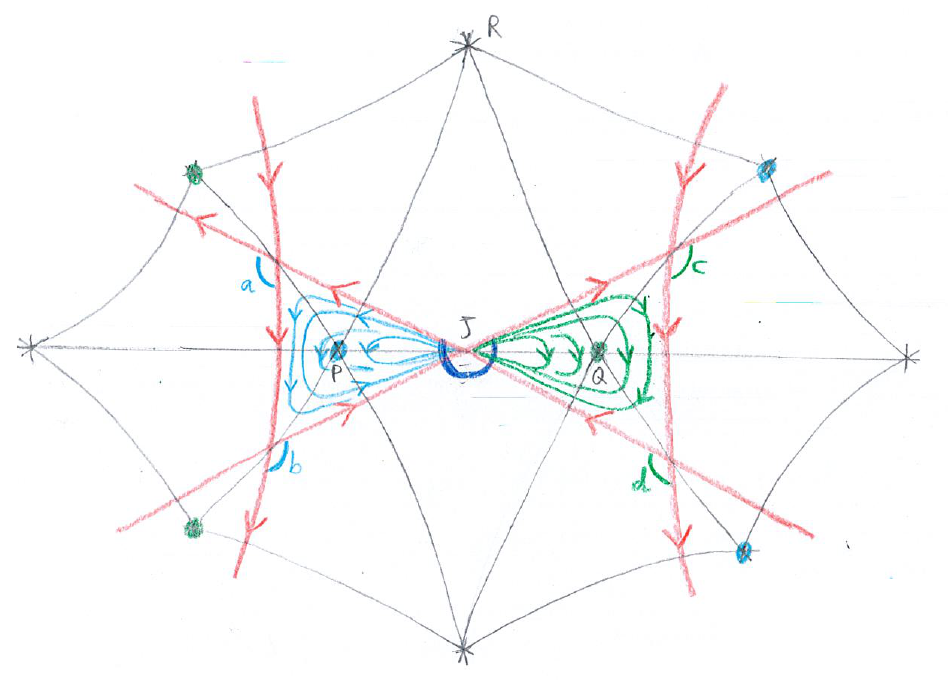}
\end{center}

Then lift~$\F_{334}$ in~$\U\Hy$ by considering the set of those unit tangent vectors positively tangent to~$\F_{334}$. 
This is a surface~$\hat S$ that is singular in the fiber of~$J$, and which is bounded by the segments of~$\vec\gamma_8$ that bound the triangles containing~$P$ and~$Q$. 

Consider the surface~$\hat S_{334}$ in~$\U\S_{334}$ obtained by projecting~$\hat S$ in~$\U\S_{334}$. 
It is also a surface that is still singular in the fiber of~$J$. 
It is transverse to~$\phi_{334}^t$, hence can be cooriented by it, and thus oriented. 
With this convention, its oriented boundary is~$-3\vec\gamma_8$. 
Indeed all sides of the triangle containing~$P$ ({\it resp}. $Q$) are identified. 

In the fiber of~$J$, the surface can be desingularized transversally to~$\phi_{334}^t$, as shown below.  
Indeed, on one part of the fiber of~$J$ (the smallest), 4 pieces arrive, 2 from the side of~$P$ ($S_a$ and $S_b$) and 2 from the side of~$Q$ ($S_c$ and $S_d$). 
One has to glue them two by two, keeping the transversality to~$\phi_{334}^t$. 
With the notation of the picture, one checks that, on the side of~$P$, the piece~$S_a$ is on top of $S_b$ (because the corresponding part of the foliation~$\F_{334}$ is in front in the trigonometric order), hence in front. 
On the side of~$Q$, the piece~$S_c$ is under the piece~$S_d$, hence in front. 
The result is a smooth surface~$S_{334}$ transverse to~$\phi_{334}$, it is bounded by~$-3\vec\gamma_8$. 

\begin{picture}(250,130)(0,0)
\put(-10,10){\includegraphics[width=.45\textwidth]{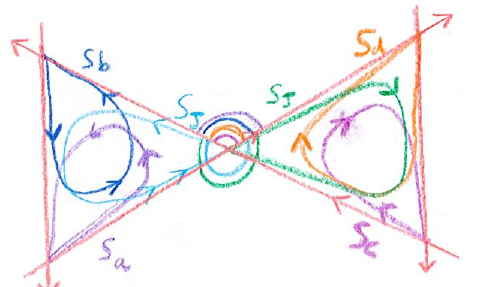}}
\put(135,0){\includegraphics[width=.22\textwidth]{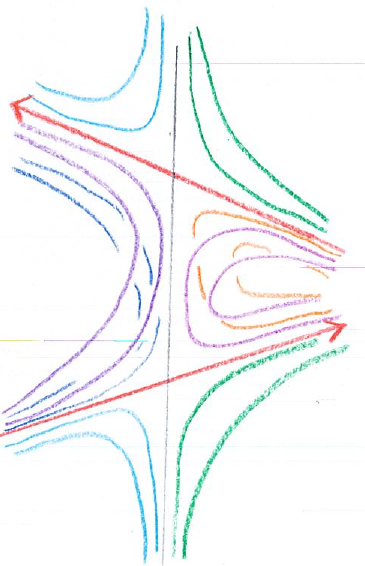}}
\put(215,12){\includegraphics[width=.28\textwidth]{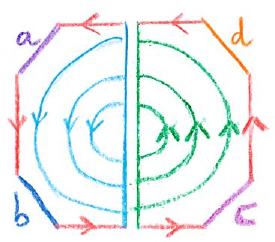}}
\end{picture}

Its genus can be computed in two ways: 
\begin{itemize}
\item either one remarks that~$S_{334}$ is made to two hexagons~$H_P, H_Q$, so that three sides of~$H_P$ are identified with three sides of~$H_Q$ following the pattern of the above picture, and so it is a torus with one boundary component;
\item or one remarks that~$S_{334}$ has one exactly boundary component, is transverse to~$\phi_{334}^t$, and winds once around its boundary in the meriodional direction. Indeed following~$\vec\gamma_8$, when encountering points $a, b, c, d$ the surfaces makes about~$-\frac12$ turn in the meridional direction and when encountering points $e, f$ it makes about $+\frac12$ turn. 
All-in-all, this makes $-1$ turn, there are just two stable and two unstable separatrices on~$S_{334}$ along its only boundary component. 
By a Poincaré-Hopf type argument,~$S_{334}$ is a one-holed torus. 
\end{itemize}

Finally one has to compute the first-return map along the geodesic flow~$\phi_{334}$ on~$S_{334}$, that is, check that there is no fixed point except the boundary orbit~$\vec\gamma_8$. 

\begin{center}
\includegraphics[height=6cm]{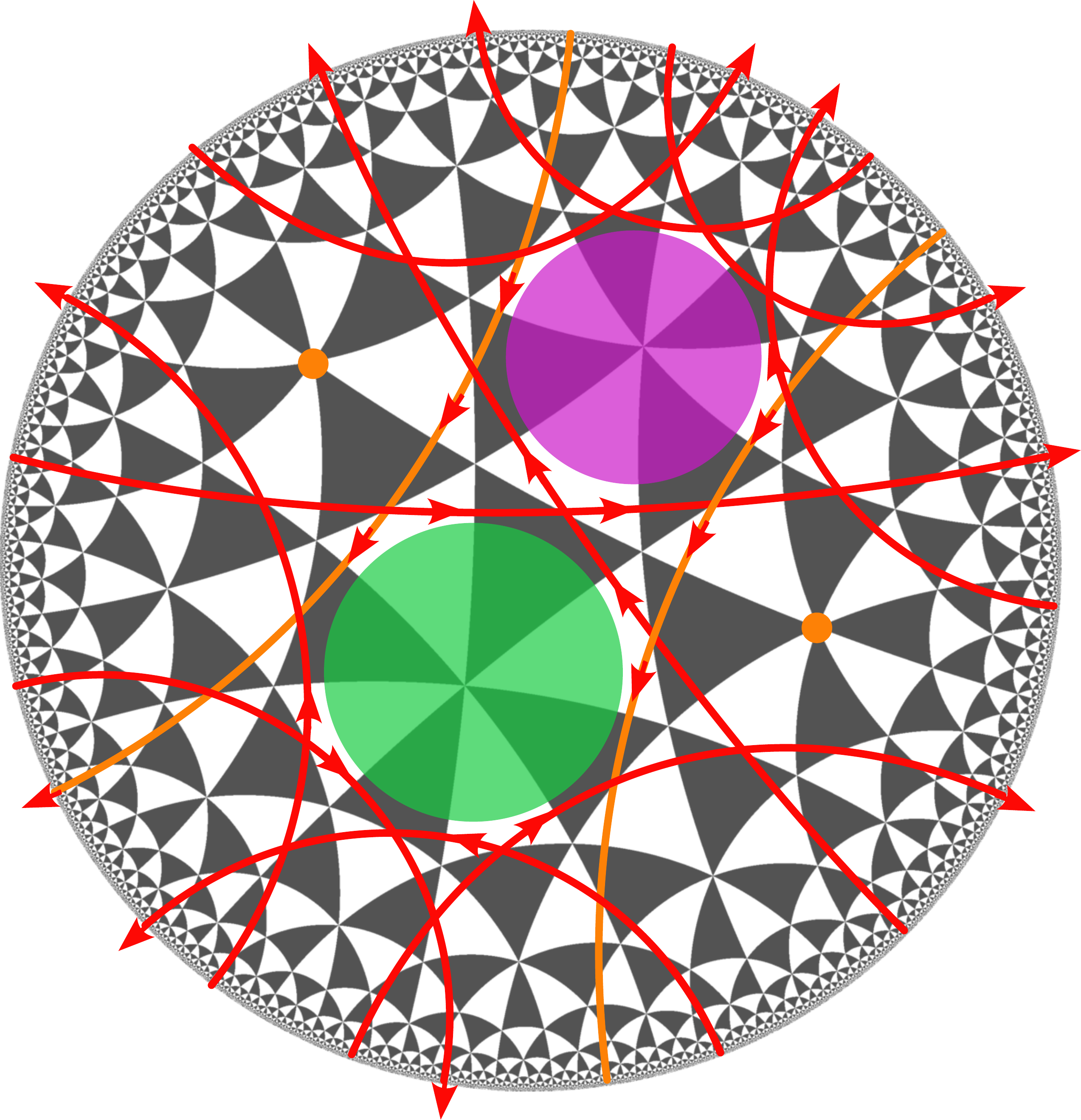}
\end{center}

As in the three previous cases, one does it by looking in~$\Hy$. 
The intersection of a periodic orbit of~$\phi_{334}^t$ with~$S_{334}$ can be read by looking at the corresponding geodesic in~$\Hy$: such a geodesic corresponds to an element~$g\in\Gamma_{334}$ that translates one lift of~$R$ into another one. 
One can decompose such a translation into two types of elementary steps: one that goes from the green octagon to the pink one on the previous picture, and one that goes from the pink one to the green one.
These two steps are different since the first one induces two intersection points with~$S_{334}$ (roughly speaking around the part $S_a=S_c$ and $S_b=S_d$ respectively), while the second one induces one intersection point only (around the part $S_J$). 
Therefore, the only way for a periodic geodesic to intersect~$S_{334}$ only once is to correspond to an element that takes the pink octagon to the green one. 
There are only 4 such elements, 2 of which are rotations around orange points, and 2 that correspond to~$\gamma$.
Thus~$f_{\bar S_{334}}$ has one fixed point only, which corresponds to the boundary component. 
Therefore it is conjugated to the cat-bat map. 


\section{The $344$-case}

This case is very similar to the previous one. 

In~$\Hy$ tiled by copies of a triangle~$PQR$ with angles $\pi/3, \pi/4, \pi/4$, call $K$ the midpoint of the segment~$QR$. 
On the orbifold~$\S_{344}$, consider the oriented geodesic~$\gamma_8$ that has one self-intersection at~$K$ and that winds once positively around~$Q$ and once negatively around~$R$. 
In~$\Hy$, the curve~$\gamma_8$ and its images by~$\Gamma_{344}$ tile~$\Hy$ into hexagons (containing the orbits under the action~$\Gamma_{344}$ of $P$), and two types of squares containing the orbits of $Q$ and $R$ respectively. 
Lift~$\gamma_8$ in~$\U\S_{344}$ into an periodic orbit~$\vec\gamma_8$ of the the geodesic flow~$\phi_{344}$. 

\begin{center}
\includegraphics[height=6cm]{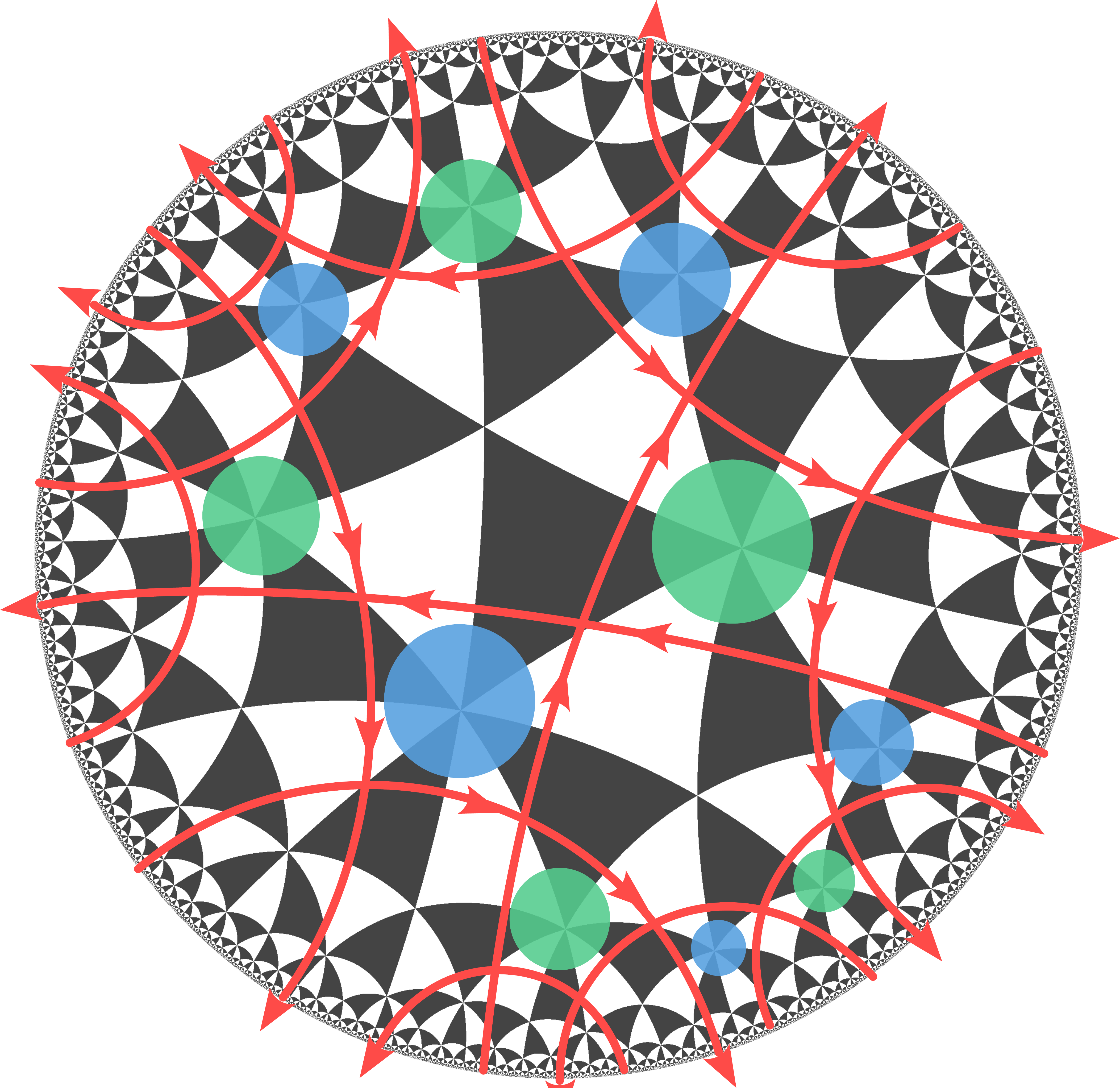}
\includegraphics[width=.45\textwidth]{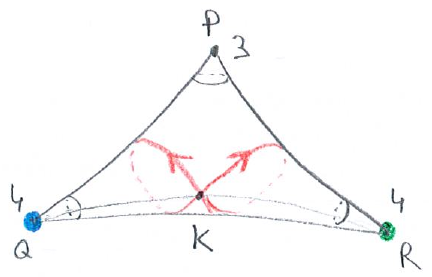}
\end{center}

In $\Hy$ consider the interiors of the two squared delimited by~$\gamma_8$ with one vertex at $K$.
These have an oriented foliation~$\F_{344}$ by convex curves tangent to~$\gamma_8$ at~$K$.

\begin{center}
\includegraphics[height=6cm]{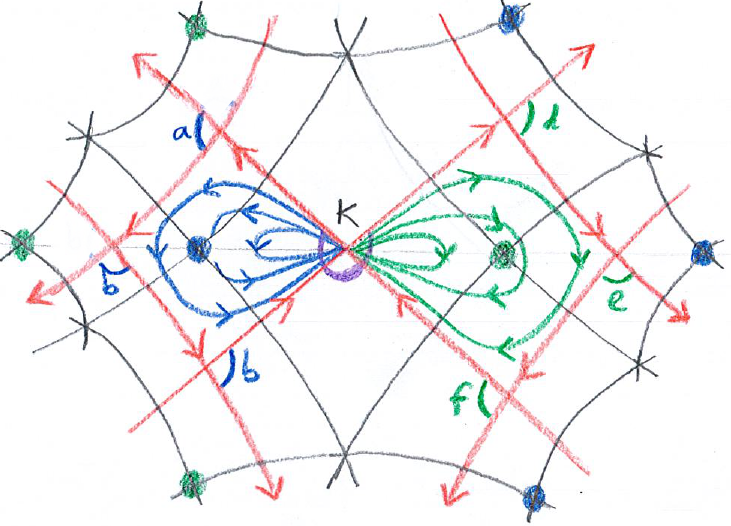}
\end{center}

Then lift~$\F_{344}$ in~$\U\Hy$ by considering the set of those unit tangent vectors positively tangent to~$\F_{344}$. 
This is a surface~$\hat S$ that is singular in the fiber of~$K$, and which is bounded by the segments of~$\vec\gamma_8$ that bound the squares containing~$Q$ and~$R$. 

Consider the surface~$\hat S_{344}$ in~$\U\S_{344}$ obtained by projecting~$\hat S$ in~$\U\S_{344}$. 
It is also a surface that is still singular in the fiber of~$K$. 
Note that its oriented boundary is~$-4\vec\gamma_8$. 
Indeed all sides of the square containing~$P$ are identified.

In the fiber of~$K$, the surface can be desingularized transversally to~$\phi_{344}^t$, as in the previous case:  
on the smallest part of the fiber of~$K$, 6 pieces arrive, 3 from the side of~$Q$ ($S_a, S_b$ and $S_c$) and 3 from the side of~$R$ ($S_d, S_e$ and $S_f$). 
One has to glue them two by two, keeping the transversality to~$\phi_{344}^t$. 
With the notation of the picture, one checks that, on the side of~$Q$, the piece~$S_a$ is on top, $S_b$ in the middle, and $S_b$ on the bottom (because or the trigonometric order of the corresponding parts of the foliation~$\F_{344}$). 
On the side of~$R$, the piece~$S_d$ is in the bottom, hence in front, $S_e$ in the middle, and $S_f$ on top, hence in the back. 
Thus one glues~$S_a$ with~$S_d$, $S_b$ with~$S_e$, and $S_c$ with~$S_f$. 
The result is a smooth surface~$S_{344}$ transverse to~$\phi_{344}^t$, it is bounded by~$-4\vec\gamma_8$. 

\begin{picture}(250,130)(0,0)
\put(-10,15){\includegraphics[height=3.2cm]{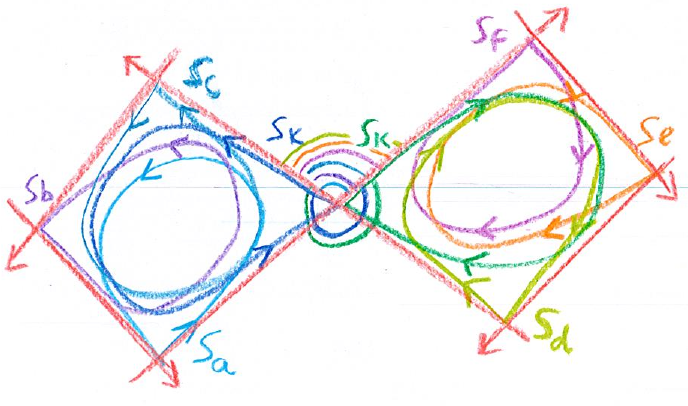}}
\put(150,0){\includegraphics[height=4.5cm]{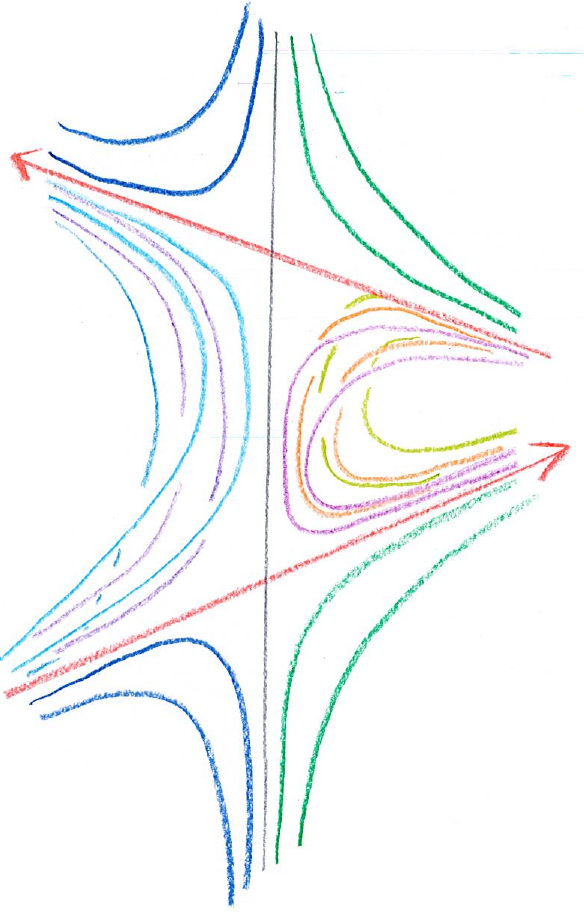}}
\put(230,20){\includegraphics[height=3cm]{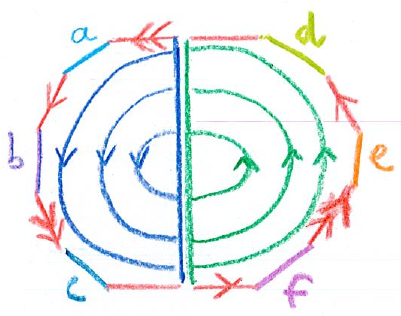}}
\end{picture}

As before, its genus can be computed in two ways: 
\begin{itemize}
\item either one remarks that~$S_{344}$ is made to two octagons~$H_Q, H_R$, so that four sides of~$H_Q$ are identified with four sides of~$H_R$ following the pattern of the above picture, and so it is a torus with two boundary components;
\item or one remarks that~$S_{344}$ has two boundary components, is transverse to~$\phi_{344}^t$, and winds twice around its boundary in the meriodional direction. 
Indeed following~$\vec\gamma_8$, when encountering points $a, b, c, d, e, f$ the surfaces makes about~$-\frac12$ turn in the meridional direction and when encountering points $g, h$ it makes about $+\frac12$ turn. 
All-in-all, this makes $-2$ turn, there are just four stable and four unstable separatrices on~$S_{344}$ along its two boundary components. 
By a Poincaré-Hopf type argument,~$S_{344}$ is a two-holed torus. 
\end{itemize}

Finally one has to compute the first-return map along the geodesic flow~$\phi_{344}^t$ on~$S_{344}$, that is, check that there is one fixed in the interior, since the boundary orbit~$\vec\gamma_8$ correspond to an orbit of period~2. 

\begin{center}
\includegraphics[height=7cm]{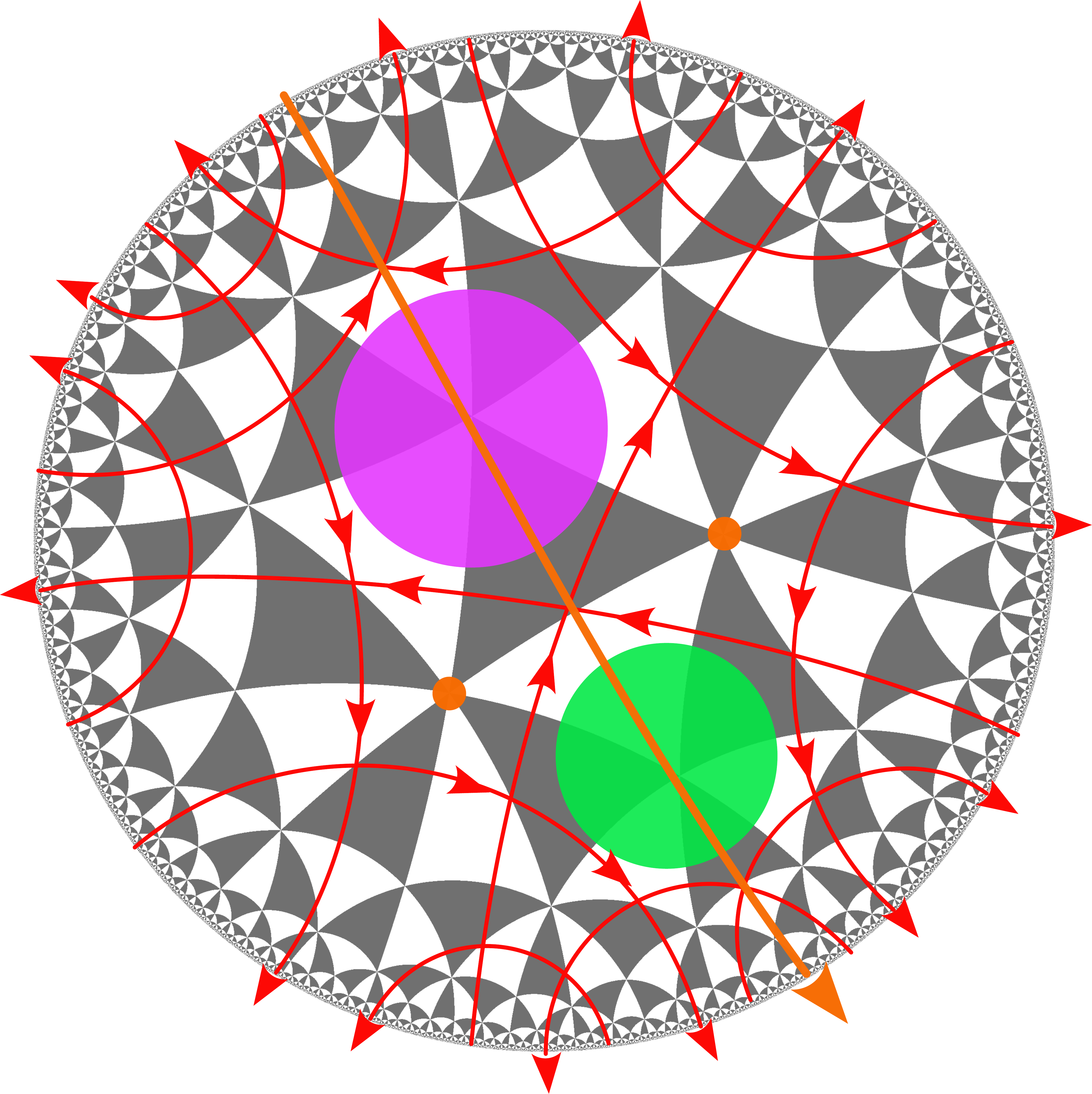}
\end{center}

As before, one does it by looking in~$\Hy$. 
The intersection of a periodic orbit of~$\phi_{344}^t$ with~$S_{344}$ can be read by looking at the corresponding geodesic in~$\Hy$: such a geodesic corresponds to an element~$g\in\Gamma_{344}$ that translates one lift of~$P$ into another one. 
One can decompose such a translation into two types of elementary steps: one that goes from the green hexagon to the pink one on the previous picture, and one that goes from the pink one to the green one.
These two steps are different since the first one induces three intersection points with~$S_{344}$ (roughly speaking around the parts $S_a=S_d$, $S_b=S_e$ and $S_c=S_f$ respectively), while the second one induces one intersection point only (around the part $S_K$). 
Therefore, the only way for a periodic geodesic to intersect~$S_{344}$ only once is to correspond to an element that takes the pink hexagon to the green one. 
There are only 3 such elements, 2 of which are rotations around orange points, and 1 that corresponds to a periodic geodesic (actually one lift of the height of~$P$ in the triangle~$PQR$, which is also the curve~$b$ in the $246$-case).
Thus~$f_{\bar S_{344}}$ has one fixed point only. 
Therefore it is conjugated to the cat-bat map. 


\bibliographystyle{siam}

\end{document}